%% file: main.tex
\documentclass[journal]{IEEEtran}
\usepackage[bookmarks=false]{hyperref}
\hypersetup{
    colorlinks=true,
    linkcolor=black,
    filecolor=black,
    urlcolor=black,
}
\usepackage[utf8]{inputenc}
\usepackage[T1]{fontenc}
\usepackage{algorithmic}
\usepackage{algorithm} 
\usepackage{nomencl}
\usepackage{etoolbox}
\usepackage[dvipsnames]{xcolor}
\usepackage{wrapfig}
\usepackage{graphicx}
\graphicspath{ {./images/} }

\newcommand{\black}[1]{\textcolor{black}{#1}}

\usepackage{tikz}
\usetikzlibrary{matrix}

\makeatletter
\newdimen\multi@col@width
\newdimen\multi@col@margin
\newcount\multi@col@count
\tikzset{
  multicol/.code={%
    \global\multi@col@count=#1\relax
    \global\let\orig@pgfmatrixendcode=\pgfmatrixendcode
    \global\let\orig@pgfmatrixemptycode=\pgfmatrixemptycode
    \def\pgfmatrixendcode##1{\orig@pgfmatrixendcode%
      ##1%
      \pgfutil@tempdima=\pgf@picmaxx
      \global\multi@col@margin=\pgf@picminx
      \advance\pgfutil@tempdima by -\pgf@picminx
      \divide\pgfutil@tempdima by #1\relax
      \global\multi@col@width=\pgfutil@tempdima
      \pgf@picmaxx=.5\multi@col@width
      \pgf@picminx=-.5\multi@col@width
      \global\pgf@picmaxx=\pgf@picmaxx
      \global\pgf@picminx=\pgf@picminx
      \gdef\multi@adjust@position{%
        \setbox\pgf@matrix@cell=\hbox\bgroup
        \hfil\hskip-\multi@col@margin
        \hfil\hskip-.5\multi@col@width
        \box\pgf@matrix@cell
        \egroup
      }%
      \gdef\multi@temp{\aftergroup\multi@adjust@position}%
      \aftergroup\multi@temp
    }
    \gdef\pgfmatrixemptycode{%
      \orig@pgfmatrixemptycode
      \global\advance\multi@col@count by -1\relax
      \global\pgf@picmaxx=.5\multi@col@width
      \global\pgf@picminx=-.5\multi@col@width
      \ifnum\multi@col@count=1\relax
       \global\let\pgfmatrixemptycode=\orig@pgfmatrixemptycode
      \fi
    }
  }
}
\makeatother
\usepackage{blindtext}
\usepackage{multicol}
\ifCLASSINFOpdf
  \graphicspath{{../pdf/}{../jpeg/}}
  \DeclareGraphicsExtensions{.pdf,.jpeg,.png}
\else

\fi

\input{comments2}

\begin{document}
\bstctlcite{IEEEexample:BSTcontrol}
\title{Unit Commitment Predictor With 
a Performance Guarantee: A Support Vector Machine Classifier}

\author{Farzaneh Pourahmadi,~\IEEEmembership{Member,~IEEE} and Jalal Kazempour,~\IEEEmembership{Senior~Member,~IEEE}

\thanks{
 The work of F. Pourahmadi  was supported in part by the Independent Research Fund Denmark (DFF) under Grant 2102-00122B. \\
F. Pourahmadi and J. Kazempour are with the Department of Wind and Energy Systems, Technical University of Denmark, Kgs. Lyngby
2800, Denmark (e-mails: \{farpour, jalal\}@dtu.dk).\
}
\vspace{-8mm}
}


\maketitle

\begin{abstract}
The system operators usually need to solve large-scale unit commitment problems within limited time frame for computation. This paper provides a pragmatic solution, showing how by learning and predicting the on/off commitment decisions of conventional units, there is a potential for system operators to warm start their solver and speed up their computation  significantly. For the prediction, we train linear and kernelized support vector machine classifiers, providing an out-of-sample performance guarantee if properly regularized, converting to distributionally robust classifiers. For the unit commitment problem, we solve a mixed-integer second-order cone problem. Our results based on the IEEE 6- and 118-bus test systems show that the kernelized SVM with proper regularization outperforms other classifiers, reducing the computational time by a factor of 1.7. In addition, if there is a tight computational limit, while the unit commitment problem without warm start is far away from the optimal solution, its warmly-started version can be solved to (near) optimality within the time limit. 

\end{abstract}


%
\IEEEpeerreviewmaketitle

\section{Introduction}
Power system operators solve the unit commitment (UC) problem on a daily basis in a forward stage, usually a day in advance, to determine the on/off commitment of conventional generating units \cite{anjos}.
They may also update the solution whenever new information, e.g., updated demand or renewable production forecasts, is available \cite{hobbs}. Despite all recent progress in developing advanced mixed-integer optimization solvers, solving the UC problem for large systems in practice can still be a computationally difficult task.
It is even getting more complex with more stochastic renewable units integrated into the system, as the operator may need to reschedule the commitment of conventional units more often than before, requiring to (re-)solve the UC problem at faster paces, likely resulting in a solution with an unsatisfactory optimality gap due to limited time window available for computation tasks. 

For example, Nord Pool collects bids every day until noon and should disseminate market outcomes not later than $2$pm, giving a maximum time window of two hours for computation tasks. Although the European markets such as Nord Pool do not solve a UC problem as it is common in the U.S. market, they still solve  a mixed-integer optimization problem due to the presence of block orders. Here, we focus on the UC problem, although our discussion on the computational complexity is also valid to the European markets.

The system and market operators  usually implement various simplifications and approximations for grid modeling, whose implications have been extensively studied over the past 50 years \cite{UC_AC1,tejada,linear,SOCP,SDP,atakan,morales,rajan,ostrowski}. However,  existing approaches could  still be computationally challenging for real large-scale systems.

\color{black}
\vspace{-1mm}
\subsection{Literature review}
\vspace{-1mm}
Several studies in the literature propose various solutions to reduce the computational complexity of the deterministic UC problem expressed as a mixed-integer linear program. Techniques employed include formulation tightening \cite{ostrowski,tightening1,tightening2,Kai}, decomposition techniques \cite{decomposition}, and constraint screening \cite{costscreening}. 
However, they overlook the fact that the same UC problem with (slightly) different input data is being solved regularly. 

A key that may give the system and market operators an advantage is that they solve the UC problem every day, and therefore they usually have access to an extensive historical database of UC solutions under various operational circumstances. This may enable the operators to speed up the computation process by \textit{learning} from previous solutions. With the recent advances in the field of machine learning, it is promising to learn from previous data to ease the computational burden.
Exploiting machine learning techniques to accelerate the rate at which mixed-integer convex optimization problems are solved is a relatively new research field \cite{Cheol}, particularly in the context of power systems. We have identified three strands in literature that use machine learning techniques for power system optimization, and more specifically for the UC problem. 

The first strand aims to develop so-called surrogate or optimization proxies, for which an estimator function is used to learn a direct mapping between contextual information, i.e., \textit{features}, such as the forecast of demand and renewable production, and the optimal operational schedules \cite{D1,Pascal_PSCC}. In the context of the UC problem, \cite{Proxy_UC1} and \cite{Proxy_UC2} develop a proxy for the UC solution by learning a direct mapping between features and optimal unit commitments.  

The second strand develops so-called indirect models, wherein an estimator function maps features to some information to reduce the optimization problem dimension. This information could be the prediction of active  constraints and/or the value of some (binary) variables. This reduces the complexity of the original optimization problem in terms of the number of constraints and/or binary variables  by eliminating inactive constraints and/or setting binary variables to predicted values \cite{ML1,roaldimplied,ID1,ID2}. In the context of the UC problem, \cite{Juanmi,UC_ActiveConstraint1,UC_ActiveConstraint2} learn the optimal value of unit commitments and active constraints to reduce the dimension of the optimization problem formulated as mixed-integer linear program. 

The main drawback of the first two strands is that they may fail to correctly predict the optimal solution. One of the main causes of this shortcoming  is that they are unable to use pre-existing mathematical form of the optimization problem as they completely replace existing optimization models by machine learning proxies. Although in the second strand,  a (reduced) optimization problem is solved to find a solution, a feasibility and optimality guarantee might be still missing.

The third strand develops an estimator function, predicting the value of some variables  as a starting point called a \textit{warm start} to speed up the solution process of non-linear or mixed-integer  problems \cite{xavier,D2,thomas}. The main advantage of this technique is that, depending on the type of the underlying optimization problem, the feasibility and likely optimality of the solution can still  be guaranteed, while the problem is being solved faster. Nevertheless, the disadvantage is that there might be cases under which a warm start brings insignificant computational benefit or even increases the computational time in an extreme case. Among others, this may happen if the warm start suggests an initial solution which is infeasible or far from the optimal point. Therefore, reducing the solution time  depends on the accuracy of warm start. Furthermore, the learning-based data-driven methodologies in the context of power systems overlook how to build an efficient training dataset and generalized learning in the case of unseen data with a performance guarantee. 
\color{black}

\subsection{Our research questions, contributions, and outline}
\vspace{-1mm}
This paper falls into the third strand which is based on warm starting. Our goal is to predict the value of binary variables, i.e., the commitment of conventional units (on/off), and use them for warm starting the integer program. 
We answer the following three questions: How can a training dataset be  built from historical UC solutions to effectively improve the warm start performance, while providing a performance guarantee? How to enhance the generalization and adaptability of learning in the case of unseen data? Lastly, we answer what the potential cost saving is through improving the computational efficiency.

To address these questions, we develop a \textit{UC predictor}, which is a tool, comprising of data collection, learning for classification, prediction and eventually decision making with a reduced computational time. To the best of our knowledge, this is the first UC predictor in the literature that  efficiently generates a training dataset and is equipped with a regularized non-linear classification approach,  resulting in a distributionally robust classifier. This predictor provides an out-of-sample performance guarantee in terms of prediction error. 

The rest of this paper is structured as follows: Section II provides an overview for the proposed three-step framework. Section III describes the data collection step. Section IV presents the learning by the classification step. Section V explains the prediction and decision making step. Section VI provides numerical results based on two case studies, including the IEEE $6$- and $118$-bus test systems.  Section VII concludes the paper. Finally, an appendix provides the UC formulation.

\vspace{2mm}
\section{Overview of the Proposed UC predictor}
\vspace{-1mm}
The proposed UC predictor is a tool with a three-step process, including (\textit{i}) data collection, (\textit{ii}) learning by classification, and (\textit{iii}) prediction and decision making. The overall structure of the proposed predictor is illustrated in Fig. 1. Each of these three steps will be explained in detail in the next three sections. However, we provide an overview in the following: 

\begin{figure} [t]
 \label{framework}
	\centering
	\includegraphics[width=1\columnwidth]{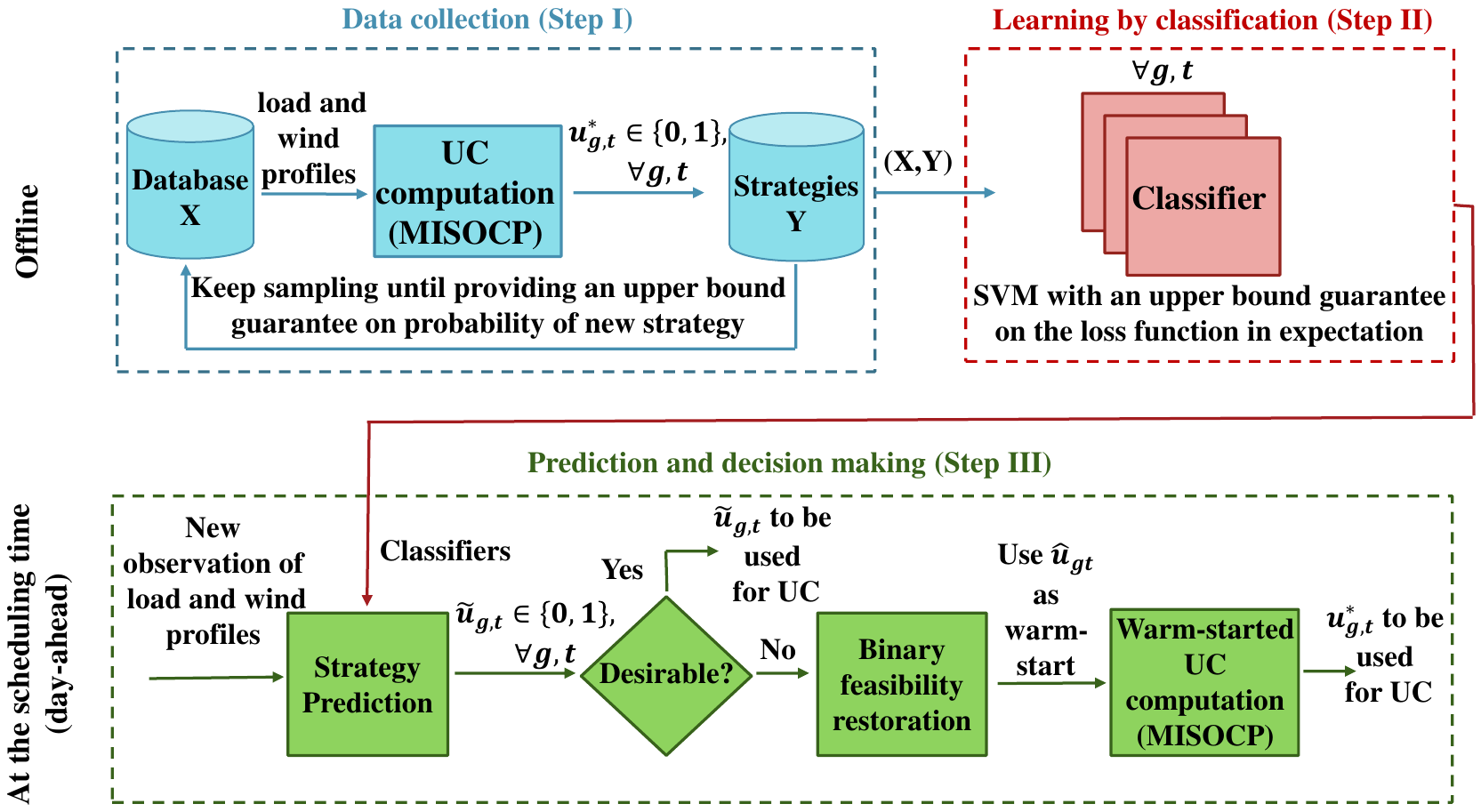}
	\caption{The overall structure of the proposed UC predictor.}
 \vspace{-2mm}
\end{figure}

The first step efficiently collects data. By data, we mean the pair of features $\mathbf{x}$ including  hourly demand and wind power forecasts over a day, and their corresponding UC solutions $\mathbf{y}$, the so-called \textit{strategies}. Each vector of features $\mathbf{x}$ results in a vector of strategies $\mathbf{y}$. Each pair ($\mathbf{x}$,$\mathbf{y}$) provides a \textit{sample}.  To accommodate the variability of demand and wind profiles, we solve the UC problem \textit{offline} multiple times, each time with different values for $\mathbf{x}$. \color{black} For the UC problem formulation, we opt for AC power flow over linearized DC power flow equations since the solutions of DC power flow are never AC feasible \cite{KyriBaker}. Aiming to include the AC power flow equations while keeping the convexity (needed for keeping the problem tractable with binaries), we use a UC formulation with the conic relaxation of AC power flow \cite{ejor}, resulting in a mixed-integer second-order cone programming (MISOCP) problem. In order to get a sufficient number of samples ($\mathbf{x}$,$\mathbf{y}$) required for the performance guarantee, we introduce 
an iterative sampling technique to construct a training dataset, ensuring an upper-bound guarantee on the probability of observing a new strategy. 
\color{black}

The second step exploits samples ($\mathbf{x}$,$\mathbf{y}$) to train a model for prediction. This model is in fact a set of \textit{binary classifiers} that provides us with a \textit{map function} in the form of $\mathbf{y} = \Phi(\mathbf{x})$. The goal is to learn what the function $\Phi(\cdot)$ is. By such a classifier, when developed, one can predict $\mathbf{y}$, i.e., the commitment of every conventional generator, when observing a new $\mathbf{x}$, without solving the UC problem. As mentioned earlier, we will use the predicted $\mathbf{y}$ to warm start the UC problem. One can hypothesize that efficient warm starting improves the branch-and-bound process and eventually accelerates the solution of an MISOCP solver. To develop each binary classifier, we use a support vector machine (SVM) technique. The reason for our choice  is its convexity and having a continuous Lipschitz loss function, enabling us to tractably reformulate the more advanced versions of the classifier. We start with a linear SVM, and then generalize the classifier by adding a regularization term. This ends up in modeling classifier as a distributionally robust optimization, providing us with a rigorous out-of-sample guarantee. We further improve the classifier by developing a kernelized SVM, aiming  to avoid under-fitting and being capable of learning in high-dimensional feature spaces.



The third step uses predicted strategies $\mathbf{y}$. If the underlying strategy is desirable, it is used as a decision for the commitment of conventional units --- therefore, no need to solve the UC problem. Otherwise, it is used as a warm start to solve the UC problem. By desirability, we refer to the solution feasibility and prediction loss guarantee. \color{black} If the prediction is undesirable due to infeasibility, we first recover a feasible solution from the sampled strategies to the infeasible one using a heuristic approach, and then use it as a warm start to the original UC problem. \color{black}

To evaluate the performance of the proposed predictor, we  solve the original and warm-started UC problem while assuming there is  a computational time limit, e.g., $15$ minutes. 
If the computational time reaches to the limit, we stop the solver and report the results obtained by then, which might be a sub-optimal solution with a high optimality gap. Our thorough out-of-sample simulations show how crucial it is to improve the accuracy of warm starting. By doing so, there is potential to reduce the system cost by solving the UC problem to (near) optimality within a limited computational time frame.

\section{Step 1: Data Collection}
\vspace{-1mm}
This section explains Step 1 of Fig. 1 (blue part).  
Suppose we have access to historical forecasts of wind and demand for $T$=$24$ hours of a day, so we can build a feature vector in the form of $\mathbf{x}=[\mathbf{p}^{\rm{W}}_1,\mathbf{q}^{\rm{W}}_1,\mathbf{p}^{\rm{D}}_1,\mathbf{q}^{\rm{D}}_1,...,\mathbf{p}^{\rm{W}}_T,\mathbf{q}^{\rm{W}}_T,\mathbf{p}^{\rm{D}}_T,\mathbf{q}^{\rm{D}}_T]\in\mathbb{R}^{4NT}$, where $N$ is the number of buses. Vectors $\mathbf{p}^{\rm{W}}$, $\mathbf{q}^{\rm{W}}$,  $\mathbf{p}^{\rm{D}}$, and $\mathbf{q}^{\rm{D}}$ refer to active and reactive power production (consumption) of wind farms (demands), respectively. For a given hour $t$, we define the active power vector  $\mathbf{p}^{\rm{W}}_t=\{p^{\rm{W}}_{1t},...,p^{\rm{W}}_{Nt}\}\in \R^N$. Vectors $\mathbf{q}_t^{\rm{W}}$,  $\mathbf{p}_t^{\rm{D}}$, and $\mathbf{q}_t^{\rm{D}}$ are defined similarly. Having access to historical data of $H$ days provides us with $H$ number of vectors $\mathbf{x}$, collected within the set of vectors $\mathbf{X}= \{\mathbf{x}_1, ..., \mathbf{x}_h, ..., \mathbf{x}_H\}$.

Given every vector $\mathbf{x}$, we solve the UC problem as a MISOCP problem, whose formulation is given in the appendix. This gives us the set of strategy vectors $\mathbf{Y}= \{\mathbf{y}_1$,$ ...$,$ \mathbf{y}_h$,$ ...$,$ \mathbf{y}_H\}$. \color{black} Every strategy vector $\mathbf{y}$ is defined as $\mathbf{y}=[u^{*}_{11}$,$...$,$u^{*}_{MT}]\in\{0,1\}^{MT}$, where $M$ is the number of conventional units, and superscript $^{*}$ refers to the optimal value obtained by the state-of-the-art solvers, e.g., Gurobi. By the optimal value obtained, we refer to the value achieved by the solver under no computational time constraint. \color{black} By this, we eventually collect the set of samples $(\mathbf{X,Y})= \{(\mathbf{x}_1,\mathbf{y}_1), ..., (\mathbf{x}_h,\mathbf{y}_h), ..., (\mathbf{x}_H,\mathbf{y}_H)\}$.

To determine the minimum required number of samples $H$ to be used later for classification, we calculate the probability of encountering a novel strategy $\mathbf{y}_{H+1}$, as in the following theorem outlined in \cite{ML1}: 

\begin{theorem} 
Given the set of feature vectors $\mathbf{X}$ drawn from an unknown discrete distribution and its corresponding set of strategy vectors $\mathbf{Y}$, the probability of coming across a feature vector $\mathbf{x}_{H+1}$ that corresponds to an unobserved strategy vector $\mathbf{y}_{H+1}$ is 
%
%
\begin{align}
    \mathbb{P}(\mathbf{y}_{H+1}\notin \mathbf{Y})\leq \frac{H_1}{H} + \tau \sqrt{\frac{1}{H}ln(\frac{3}{\epsilon})}, \label{prob}
\end{align} 
  where $H_1$ is the number of unique strategies that have been observed only once, $\tau=2\sqrt{2}+\sqrt{3}$, and $1-\epsilon$ is the confidence level of holding (\ref{prob}).
\end{theorem}

\begin{algorithm} [t]
\caption{Strategy sampling} 
Given $\delta$ and $\epsilon$ \\
Initialize $H=1$ and $H_1=1$\\
\textbf{While} $\frac{H_1}{H} + \tau \sqrt{\frac{1}{H}ln(\frac{3}{\epsilon})} > \delta$

\emph{a)} Sample $\mathbf{x}_h$ and compute $\mathbf{y}_{h}$ \\
\emph{b)} Update the set of features $\mathbf{X} \leftarrow \mathbf{X} \cup \{\mathbf{x}_{h}\}$ \\
\emph{c)}  Update the set of corresponding strategies $\mathbf{Y} \leftarrow \mathbf{Y} \cup \{\mathbf{y}_{h}\}$ \\
\emph{d)} \textbf{if} $\mathbf{y}_{h} \notin \mathbf{Y}$ \textbf{then} $H_1=H_1+1$\\
\emph{e)} $H=H+1$\\
\textbf{end}\\ 
\textbf{Return} $\mathbf{X}$ and $\mathbf{Y}$
\end{algorithm}

The proof can be directly derived from Theorem $9$ of \cite{Schapire}. \color{black}
This theorem provides us with an approach to estimate the probability of observing an unseen set of commitment decisions. 
Let us delve deeper into this theorem. Consider a set of commitment solutions, denoted as $\mathbf{Y}$, with a total of $H$ solutions. Let $\mathbf{Y}_k$ denote the subset of commitment solutions occurring exactly $k$ times within $\mathbf{Y}$. We use $M_k$ to represent the probability of observing a commitment solution $k$ times. Good-Turing estimators, denoted as $G_k$, provide estimates for $M_k$ as 
\begin{equation}
    G_k=\frac{k+1}{H-k}|\mathbf{Y}_{k+1}|, 
\end{equation}
where $|\cdot|$ denotes the size of a set. In particular, the estimator $G_0$ for observing unseen commitment solutions is given by $G_0=\frac{|\mathbf{Y}_{1}|}{H}$. We denote $|\mathbf{Y}_{1}|$ as $H_1$, representing the count of unique commitment solutions or strategies. Intuitively, when there are numerous different commitment solutions but a relatively small training dataset, the estimate of $M_0$ may be 1. This scenario can occur when collecting a large number of training samples is computationally burdensome. Therefore, $G_0$ is not a good estimate of $M_0$. In order to resolve this problem, \cite{Schapire} establishes a bound on the bias of the Good-Turing estimator $G_0$ and extends it using confidence notations. It provides a tight upper bound on the special case of $M_0$ with a confidence level of $1-\epsilon$, expressed as

\begin{equation} \label{Good_Turing}
    M_0 \leq \frac{H_1}{H} + (2\sqrt{2}+\sqrt{3})\sqrt{\frac{1}{H}ln(\frac{3}{\epsilon})}.
\end{equation}  

In other words,  (\ref{Good_Turing}) is satisfied with a probability of at least $1-\epsilon$ over the choice of samples. Considering this theorem, we keep sampling until the bound on the right-hand side of (\ref{Good_Turing}) falls below a desired probability guarantee of $\delta$ with a given confidence interval of $1-\epsilon$. This strategy sampling is detailed in Algorithm 1.

\color{black}
\textit{Remark}: 
Hereafter, in order to train the classifiers, we convert $y=0$ in our samples to $y=-1$.  Therefore, $y_{gt}=-1$ implies that  unit $g$ in hour $t$ is off. In contrast, $y_{gt}=1$ still means the unit is on.

\section{Step 2: Learning by Binary Classification}
This section describes Step 2 of Fig. 1 (red part). 
Given the training dataset, i.e., the set of samples, $(\mathbf{X,Y})= \{(\mathbf{x}_1,\mathbf{y}_1),..., (\mathbf{x}_h,\mathbf{y}_h),..., (\mathbf{x}_H,\mathbf{y}_H)\}$, this step designs a set of binary classifiers to construct the function $\Phi:\mathbb{R}^{4NT}\rightarrow\{-1,1\}^{MT}$ mapping $\mathbf{X}$ to $\mathbf{Y}$. As the number of strategies required to satisfy (\ref{prob}) can increase quickly for some problems, the implementation of a multi-class classification method may impose a heavy computational burden. To overcome it, we make a simplification assumption and develop one binary classifier per each entry of $\mathbf{y}$. This means that for each conventional unit $g$ and time period $t$, we identify a distinct map function $\phi_{gt}:\mathbb{R}^{4NT}\rightarrow\{-1,1\}$. We then gather all of these map functions as $\Phi=\{\phi_{11},...,\phi_{gt},...,\phi_{MT}\}$ to predict the on/off commitment status of conventional units $g=\{1,...,M\}$ at time periods $t=\{1,...,T\}$.

As already discussed in Section II, we use a SVM technique to develop our binary classifier, providing us with the map functions $\Phi$. 
We start with a linear SVM, generalize it by adding a regularizer to avoid over-fitting, and finally extend it to a high-dimensional feature space as a non-linear classifier, resulting in a kernelized SVM. 


\vspace{-1mm}
\subsection{Linear SVM}
Recall we develop one classifier per conventional unit $g$ per hour $t$. Hereafter, we focus on developing a linear SVM for unit $g$ in hour $t$. For notational simplicity, we drop indices $g$ and $t$. 
For every sample $h=\{1,...,H\}$, the vector $\mathbf{x}_h \in \R^{4N}$ and the binary value $y_h \in \{-1,1\}$ are input data for the linear SVM classifier. The goal is to divide $H$ samples into two classes by a hyperplane, maximizing the distance between the hyperplane and the nearest point  \cite{ML,kernel2}. For that, we solve a linear  problem as
\begin{subequations} \label{Linear_SVM_1}
\begin{align} 
\min_{\zeta_h \geq 0,\mathbf{w},b}& \ \ \frac{1}{H} \sum_{h=1}^H \zeta_h \\
     \text{s.t.} \  \ & y_h (\mathbf{w}^\top \mathbf{x}_h+ b) \ge 1-\zeta_h, \ \forall h, \label{soft}\
\end{align}
\end{subequations}
where $\mathbf{w}\in \R^{4N}$ and $b\in \R$ constitute the hyperplane, and  $\zeta_h\in \R^+$ is a slack variable. 

Once the optimal values for  $\mathbf{w}^{*}$ and $b^{*}$ are obtained by solving \eqref{Linear_SVM_1}, any new feature vector $\tilde{\mathbf{x}}$ is classified using the map function 
\begin{align} 
\phi(\tilde{\mathbf{x}})=\text{sign}(\mathbf{w}^{{*}^\top}\tilde{\mathbf{x}}+b^{*}), \label{map1}
\end{align}
where sign($a$) =1 if $a>0$, and $-1$ if $a<0$. In other words, $\phi(\tilde{\mathbf{x}})$ predicts $\tilde{y}$, such that $\tilde{y}=1$ means the classifier predicts the unit will be on, and $\tilde{y}=-1$ otherwise. 
%
%
\vspace{-2mm}
\subsection{Reformulating Linear SVM}
We reformulate \eqref{Linear_SVM_1} to build up the next section accordingly. The optimal solution to the linear problem \eqref{Linear_SVM_1} would be in one of the corner points, such that in the optimal solution it holds
\begin{align} \label{zeta}
    \zeta_h=\max\Big(0,1-y_h(\mathbf{w}^\top \mathbf{x}_h+ b)\Big), \ \forall{h}.
\end{align}

Defining the hinge loss function $\mathbf{L}(z,z^{\prime})=\max \big(0,1-zz^{\prime}\big)$, and replacing $\zeta_h$ in (\ref{Linear_SVM_1}) with the right hand side of (\ref{zeta}), we get
\begin{align} \label{SVM_regularized}
&\min_{\mathbf{w},b} \  \frac{1}{H} \sum_{h=1}^H \mathbf{L}\big(y_h,(\mathbf{w}^\top \mathbf{x}_h +\! b )\big),
\end{align} 
which is an unconstrained linear optimization problem, equivalent to (\ref{Linear_SVM_1}). If we assume the feature vector $\mathbf x_h$ follows a probability distribution $\hat{\mathbb{P}}_H$, (\ref{SVM_regularized}) can be written as
\begin{align}\label{SVM_regularized_continous}
\min_{\mathbf{w},b} \ \mathbb{E}^{\hat{\mathbb{P}}_H}  \bigg[\mathbf{L}\big(y,(\mathbf{w}^\top\! \mathbf{x}+ b )\big)\bigg].
\end{align} 

Optimization (\ref{SVM_regularized_continous}) is an alternative derivation of the linear SVM minimizing the expected risk, where $\mathbb{E}^{\hat{\mathbb{P}}_H}[\cdot]$ is the expectation operator with respect to the empirical distribution $\hat{\mathbb{P}}_H$ of training dataset with $H$ samples. In the following section, we use  \eqref{SVM_regularized_continous} to develop a distributionally robust classifier.

\color{black}
\vspace{-1mm}
\subsection{Distributionally robust  SVM}

The literature often assumes that training and test datasets follow the same distribution, however this is not necessarily a valid assumption. When the distribution of training data is unknown, it is helpful to create a classifier which is robust to distributional uncertainty. One potential method is to use a Wasserstein distributionally robust optimization, yielding a distributionally robust SVM.
Built upon \eqref{SVM_regularized_continous}, the distributionally robust  SVM writes as
\begin{align} \label{DRC_linear}
\min_{\mathbf{w},b}& \ \sup_{\mathbb{Q}\in \mathbb{B}_\rho(\hat{\mathbb{P}}_H)} \mathbb{E}^\mathbb{Q} \bigg[\mathbf{L}\big(y,(\mathbf{w}^\top \mathbf{x}+ b )\big)\bigg],
\end{align} 
where $\mathbb{B}_\rho(\hat{\mathbb{P}}_H)$ denotes the set of distributions, whose distance to the empirical distribution $\hat{\mathbb{P}}_H$ is lower than or equal to the pre-defined value $\rho\in \R$. Optimization \eqref{DRC_linear} obtains the worst distribution $\mathbb{Q}$ within $\mathbb{B}_\rho(\hat{\mathbb{P}}_H)$, and optimally classifies samples against it. It has been proven in \cite{Soroush1} and \cite{Soroush2} that, under certain circumstances, \eqref{DRC_linear}  is equivalent to a \textit{regularized} version of the linear SVM \eqref{SVM_regularized_continous}. 
%
%
Accordingly, for every $\rho \ge 0$, there exists a value $\lambda\in \R$ such that the distributionally robust  SVM is analogous to the regularized  SVM. 
As a result, \eqref{DRC_linear} is equivalently rewritten as
\begin{align} \label{SVM_DRC_regularized}
\min_{\mathbf{w},\mathbf{b}}& \ \frac{1}{H}  \sum_{h=1}^H \mathbf{L}\big(y_h,(\mathbf{w}^\top \mathbf{x}_h+b )\big)+\lambda \lVert \mathbf{w} \rVert ^2_2,
\end{align} 
where  $\lVert \cdot\rVert_2$ is the $2$-norm. The second term in (\ref{SVM_DRC_regularized}) is the regularization term to avoid over-fitting. 
The resulting distributionally robust SVM offers an upper-bound guarantee on the expected prediction loss function, commonly referred to as an out-of-sample guarantee.
Let us consider a test dataset $\{(\Tilde{\mathbf{x}}_1,\Tilde{y}_1),...,(\Tilde{\mathbf{x}}_h,\Tilde{y}_h),...,(\Tilde{\mathbf{x}}_{\Tilde{H}},{\Tilde{y}}_{\Tilde{H}})\}$ of $\Tilde{H}$ samples. Assuming the distribution of the test dataset lies in $\mathbb{B}_\rho(\hat{\mathbb{P}}_H)$, it holds that
%
\begin{align} \label{predictionguarantee_linear}
 \frac{1}{\Tilde{H}} \sum_{h=1}^{\Tilde{H}} \mathbf{L}\big(\Tilde{y}_h,(\mathbf{w^{*}}^\top \Tilde{\mathbf{x}}_h+ b^{*} )\big)  \leq J(\lambda),
\end{align}
where $J(\lambda)$ is the optimal value of (\ref{SVM_DRC_regularized}). 
%
Problem (\ref{SVM_DRC_regularized}) can be efficiently solved by replacing $\lambda \lVert \mathbf{w} \rVert ^2_2 $  by an auxiliary variable $\varrho$ as
%
\begin{subequations} \label{Linear_SVM_2}
\begin{align} 
\min_{\zeta_h \geq 0,\mathbf{w},b}& \ \ \frac{1}{H} \sum_{h=1}^H \zeta_h + \varrho \\
     \text{s.t.} \  \ & y_h (\mathbf{w}^\top \mathbf{x}_h+ b) \ge 1-\zeta_h, \ \forall h \\
     & \varrho \geq \lambda \lVert \mathbf{w} \rVert ^2_2.
     \label{soft}\
\end{align}
\end{subequations}

 By doing so, \eqref{Linear_SVM_2} becomes a second-order cone programming (SOCP) problem with a linear objective function.

\color{black}
\vspace{-1mm}
\subsection{Kernelized SVM}
By solving either the linear SVM \eqref{Linear_SVM_1} or the  distributionally robust SVM (\ref{SVM_DRC_regularized}) or its equivalent regularized SVM \eqref{Linear_SVM_2}, we end up in a linear map function $\phi({\mathbf{x}})=\text{sign}(\mathbf{w}^{{*}^\top}{\mathbf{x}}+b^{*})$. Aiming to a more efficient classification with a higher degree of freedom, it is desirable to obtain a non-linear map function while preserving the convexity of optimization problems. For this purpose, we use the \textit{kernel trick} to transform data into a higher dimensional space.


The rationale is to transform the regularized SVM \eqref{Linear_SVM_2} from the Euclidean space to an infinite-dimensional Hilbert space $\mathbb{H}$, obtaining a linear map function there, and transforming it back to the Euclidean space as a non-linear function \cite{kernel2,patternreco}. It has been discussed in the literature that it is computationally desirable to transform the dual optimization of \eqref{Linear_SVM_2}. This dual optimization in the Euclidean space is 
\begin{subequations}\label{dual}
\begin{align} 
\max_{\alpha_h} & \sum_h \alpha_h - \frac{1}{2} \sum_{k,h} \alpha_k \alpha_h y_k y_h  \mathbf{x}_k\mathbf{x}_h \label{dual_obj} \\
\text{s.t.} \  & \sum_h y_h \alpha_h =0 \label{dual_con1}\\
& 0\leq \alpha_h \leq \lambda, \ \forall h, \label{dual_con2}
\end{align}
\end{subequations}
where $h,k=\{1,...,H\}$ are indices for samples. In addition, $\alpha_h$ is the dual variable corresponding to \eqref{soft}. This dual SVM is the Hilbert space \cite{vapnik} is written as 
\begin{subequations}
\begin{align} 
\max_{\alpha_h} & \sum_h \alpha_h - \frac{1}{2} \sum_{k,h} \alpha_k \alpha_h y_k y_h  \Psi(\mathbf{x}_k)\Psi(\mathbf{x}_h)  \label{dual_newspace_obj} \\
\text{s.t.}\  &{\rm{\eqref{dual_con1}-\eqref{dual_con2}}}, 
\end{align}
\end{subequations}
where the map $\Psi:\R^{4N}\rightarrow \mathbb{H}$ transforms the vectors from the Euclidean to the Hilbert space. The kernel trick states that there are functions $f(\mathbf{x}_k,\mathbf{x}_h)$ in the Euclidean space which are equivalent to the dot product $\Psi(\mathbf{x}_k)\Psi(\mathbf{x}_h)$ in the Hilbert space. One popular function satisfying such a property is the Gaussian kernel function \cite{kernel1,Keerthi}, defined as $f(\mathbf{x}_k,\mathbf{x}_h)=e^{-\gamma \lVert \mathbf{x}_k-\mathbf{x}_h \rVert ^2_2}$,
where the parameter $\gamma$ controls the flexibility degree of the hyperplane (tuned by cross-validation). Now, the kernelized SVM (dual form) in the Euclidean space is written as
\begin{subequations} \label{dual}
\begin{align} \label{dual_kernel}
\max_{\alpha_h} & \sum_h \alpha_h - \frac{1}{2} \sum_{k,h} \alpha_k \alpha_h y_k y_h  f(\mathbf{x}_k,\mathbf{x}_h) \\
\text{s.t.}\  &{\rm{\eqref{dual_con1}-\eqref{dual_con2}}}. 
\end{align}
\end{subequations}

Now, we convert \eqref{dual} back to its primal form: 
\begin{subequations} \label{dual_kernel_primal}
\begin{align} 
    \min_{\zeta_h\geq 0, \boldsymbol\beta} & \ \ \frac{1}{H} \sum_{h=1}^H \zeta_h +\lambda \lVert {\boldsymbol f^{\frac{1}{2}}\boldsymbol \beta} \rVert ^2_2\\
     \text{s.t.} \  & \sum_{k=1}^H y_k f(\mathbf{x}_k,\mathbf{x}_h) \beta_k  \ge 1-\zeta_h, \ \forall h, \label{SOC_con1}\
\end{align}
\end{subequations}
where $\boldsymbol\beta\in \R^{H}$ and $\boldsymbol f\in \R^{H \times H}$. \color{black} Note that $\boldsymbol f\in \R^{H \times H}$ is a positive semidefinite matrix that can be written as $\boldsymbol f=\boldsymbol f^{\frac{1}{2}} \boldsymbol f^{\frac{1}{2}}$. We use the Cholesky decomposition to find  $\boldsymbol f^{\frac{1}{2}}\in \R^{H \times H}$. \color{black} Similar to \eqref{Linear_SVM_2}, optimization \eqref{dual_kernel_primal} is a QP problem, however one can convert it to a SOCP problem by moving the quadratic regularizer from the objective function to constraints. 

Once \eqref{dual_kernel_primal} is solved and the optimal value for  $\boldsymbol\beta^{*}$ is obtained, any new feature vector $\tilde{\mathbf{x}}$ is classified using the map function 
\begin{align} 
\phi(\tilde{\mathbf{x}})=\text{sign}\Big(\sum_{h=1}^H f(\tilde{\mathbf{x}},\mathbf{x}_h)\beta^{*}_h\Big), \label{map2}
\end{align}
where the interpretation is the same as that of \eqref{map1}.

\color{black}
\vspace{-1mm}
\subsection{Distributionally Robust Kernelized SVM}

Similar to the discussion for \eqref{DRC_linear}, it has been proven in \cite{Soroush1}  and \cite{Soroush2} that the kernelized SVM  
with the regularization term $\lambda \lVert {\boldsymbol{f}^{\frac{1}{2}}\boldsymbol \beta} \rVert ^2_2$, under certain conditions, is equivalent to a distributionally robust kernelized classifier as
\begin{align} \label{DRC_nonlinear}
\min_{\boldsymbol\beta}& \sup_{\mathbb{Q}\in \mathbb{B}_\rho(\hat{\mathbb{P}}_H)} \mathbb{E}^\mathbb{Q} \Big[\mathbf{L}\Big(y,\sum_{h=1}^H   f(\mathbf{x},\mathbf{x}_h) \beta_h\Big)\Big].
\end{align}

Similar to \eqref{predictionguarantee_linear}, it holds that
\begin{align}\label{predictionguarantee_kernel}
 \frac{1}{\Tilde{H}} \sum_{k=1}^{\Tilde{H}} \mathbf{L}\Big(\Tilde y_k,\sum_{h=1}^H   f(\Tilde{\mathbf{x}}_k,\mathbf{x}_h) \beta^{*}_h\Big) \leq J'(\lambda),
\end{align}
where $J'(\lambda)$ is the optimal value of \eqref{dual_kernel_primal}.
\color{black}

A summary of various binary classifications is given in Table \ref{classfiers}.

 \begin{table}[b]
	\centering
\caption{Summary of various binary classifications}
\label{classfiers}
\begin{tabular}{c c c}
\hline \hline
{ Type of} &   {Without} &   {With regularization} \\
{ Classification} &   {regularization} &   {(distributionally robust)} \\\hline
 {Linear SVM}    & \eqref{Linear_SVM_1} or \eqref{Linear_SVM_2} with $\lambda=0$  & \eqref{Linear_SVM_2}\\
 {Kernelized SVM} &  \eqref{dual_kernel_primal} with $\lambda=0$  &  \eqref{dual_kernel_primal}  \\  \hline \hline
\end{tabular}
\end{table}

\section{Step 3: Prediction and Decision Making}
This section explains Step 3 of Fig. 1 (green part). The map functions \eqref{map1} for the linear SVM and \eqref{map2} for the kernelized SVM enable us to predict the commitment decision $\tilde{\mathbf{y}}$ of conventional units, having features $\tilde{\mathbf{x}}$ as the input data. 
%
To guarantee the prediction performance, we present a confidence level for inaccuracy of predicted $\tilde{\mathbf{y}}$ based on the concept of distributionally robust classification. This guarantee provides an upper bound on the prediction loss in expectation, which is already mentioned in (\ref{predictionguarantee_linear}) and (\ref{predictionguarantee_kernel}). Note that we here consider the hinge loss function to measure the accuracy of prediction.

\begin{figure}[t]
	\centering
	\includegraphics[width=0.95\columnwidth]{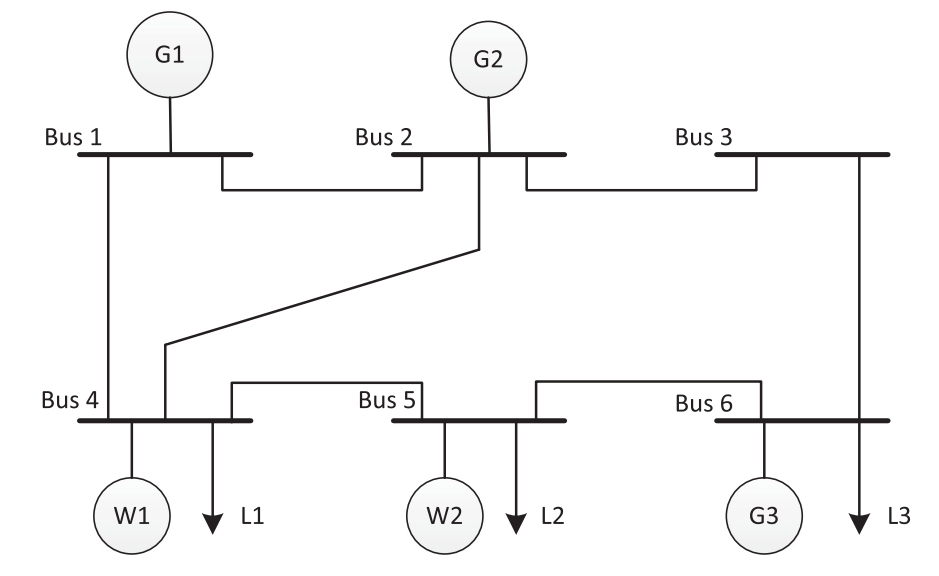}
	\caption{The IEEE $6$-bus test system with three 
conventional units $\rm{G}1$-$\rm{G}3$, two wind farms $\rm{W}1$-$\rm{W}2$, and three loads $\rm{L}1$-$\rm{L}3$.}
   \label{6bus}
\end{figure}

\begin{figure*}[t] 
  \centering
\begin{tikzpicture}
\matrix (a)[row sep=0mm, column sep=0mm,  matrix of nodes] at (0,0) {
       \textcolor{Maroon}{\textbf{Training}} & \textcolor{Maroon}{\textbf{Training}} & \textcolor{Blue}{\textbf{Testing}}& \textcolor{Blue}{\textbf{Testing}}\\
        \textcolor{Maroon}{Linear SVM} & \textcolor{Maroon}{Kernelized SVM} &\textcolor{Blue}{Linear SVM}&\textcolor{Blue}{Kernelized SVM}\\
        \includegraphics[width=4.1cm]{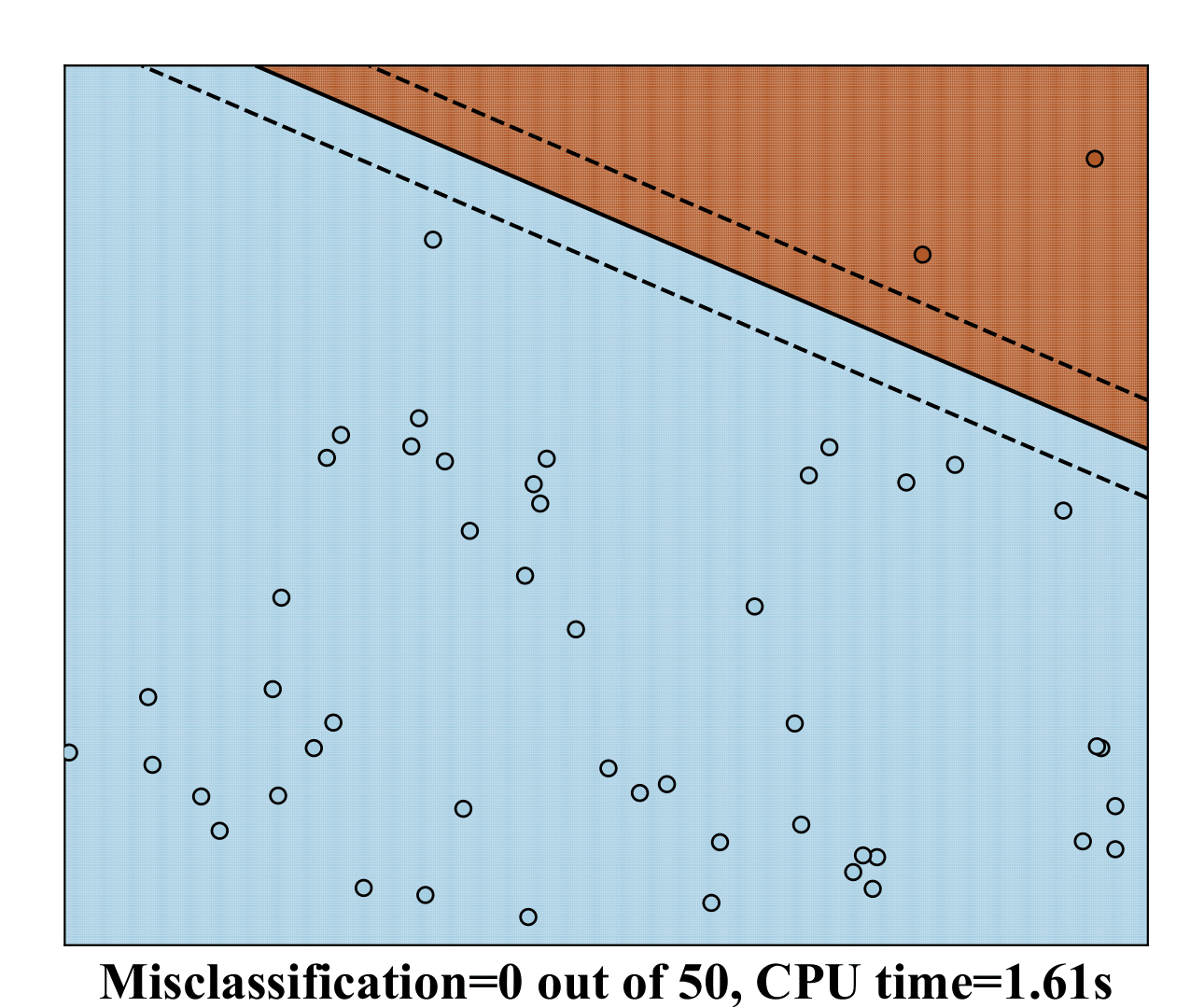} &
       \includegraphics[width=4.1cm]{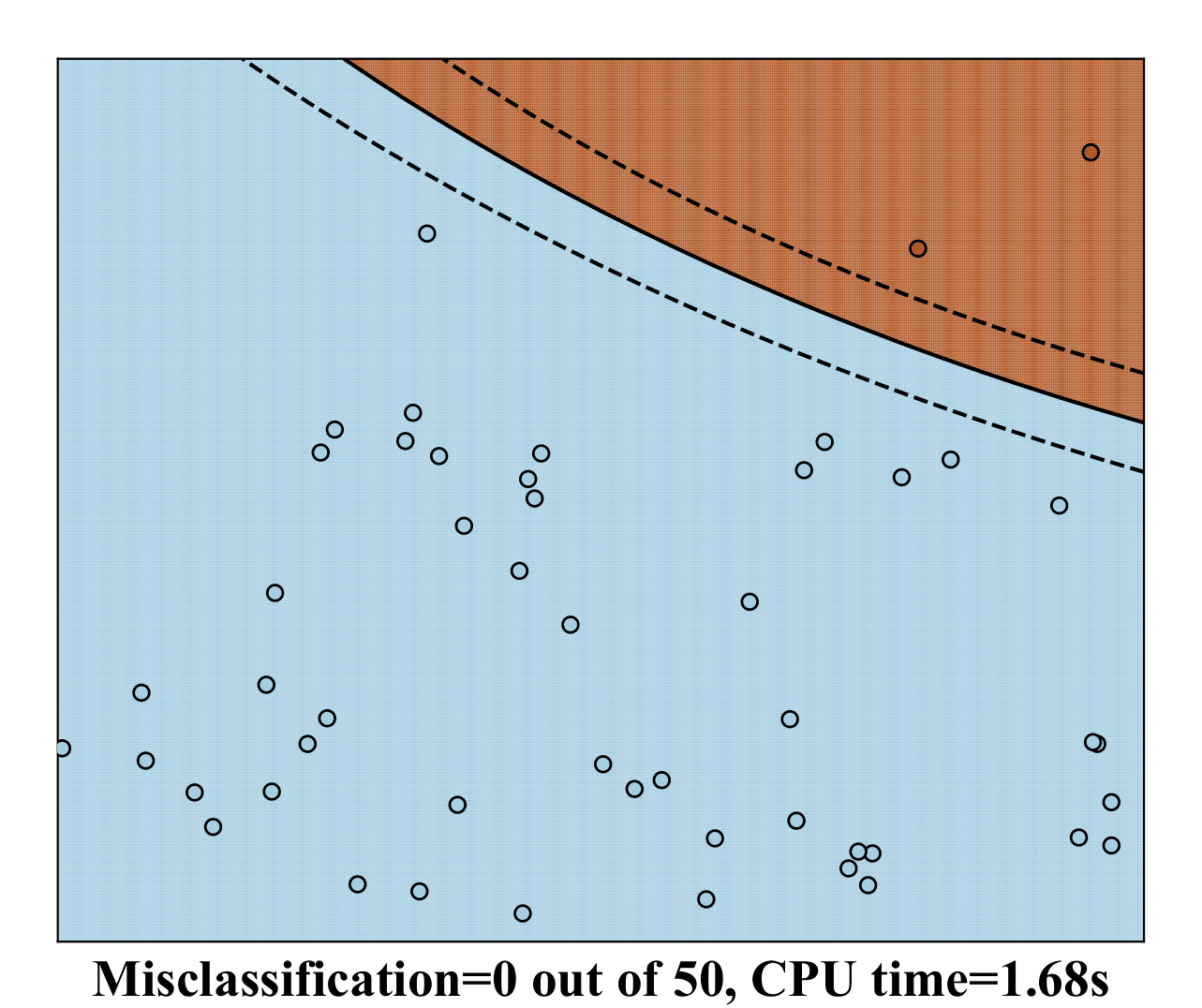} &
       \includegraphics[width=4.1cm]{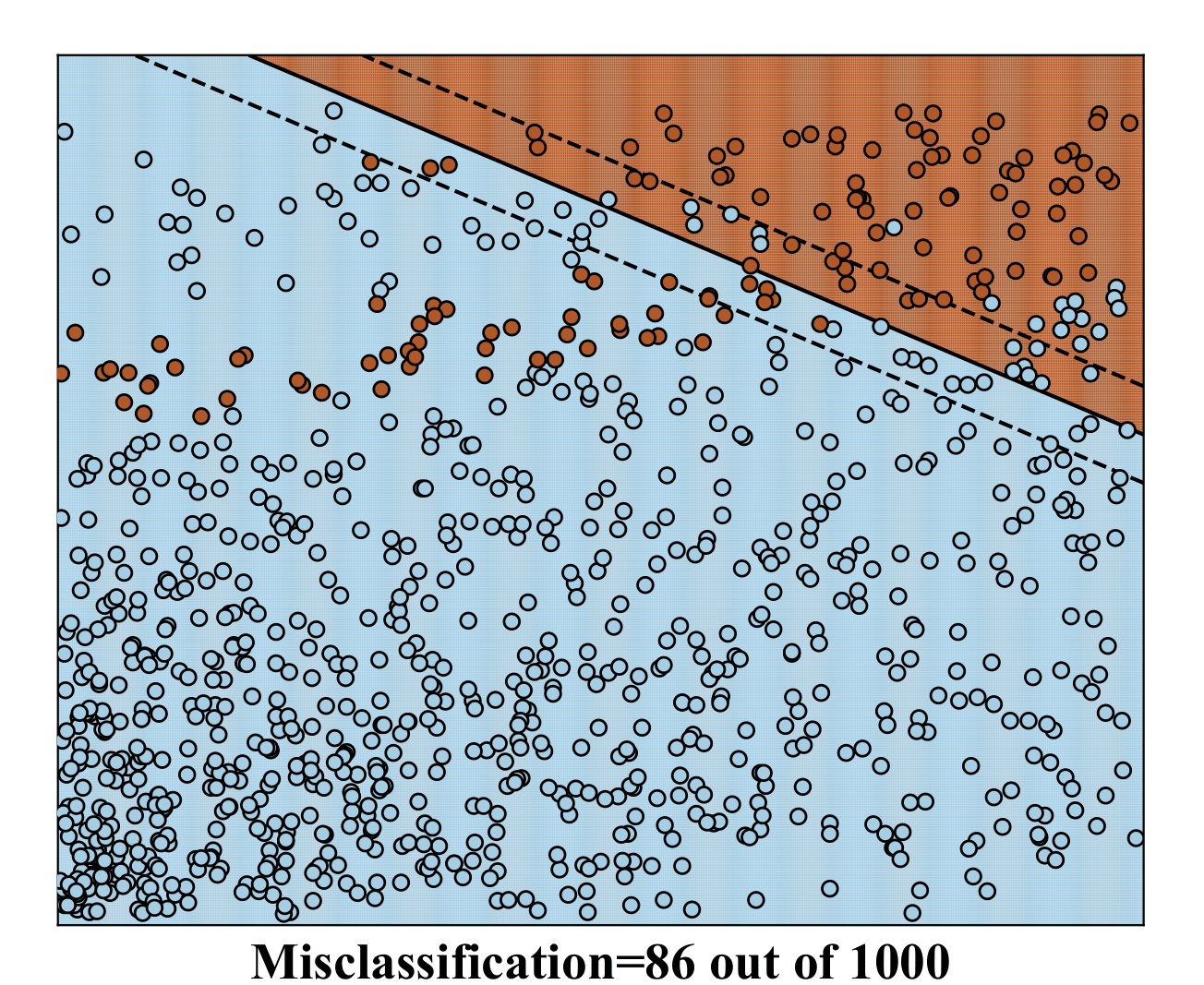} &
       \includegraphics[width=4.1cm]{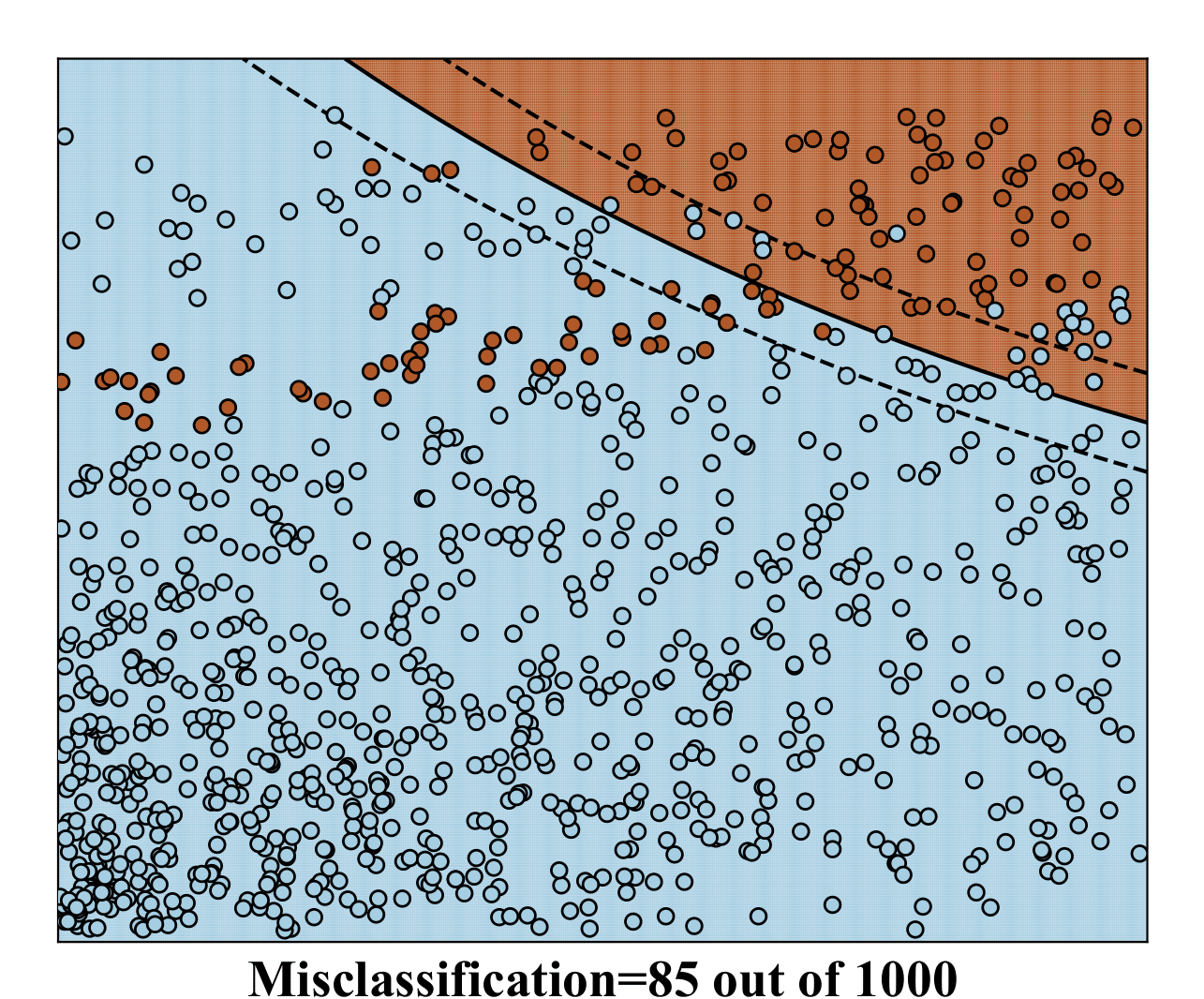} \\
        \includegraphics[width=4.1cm]{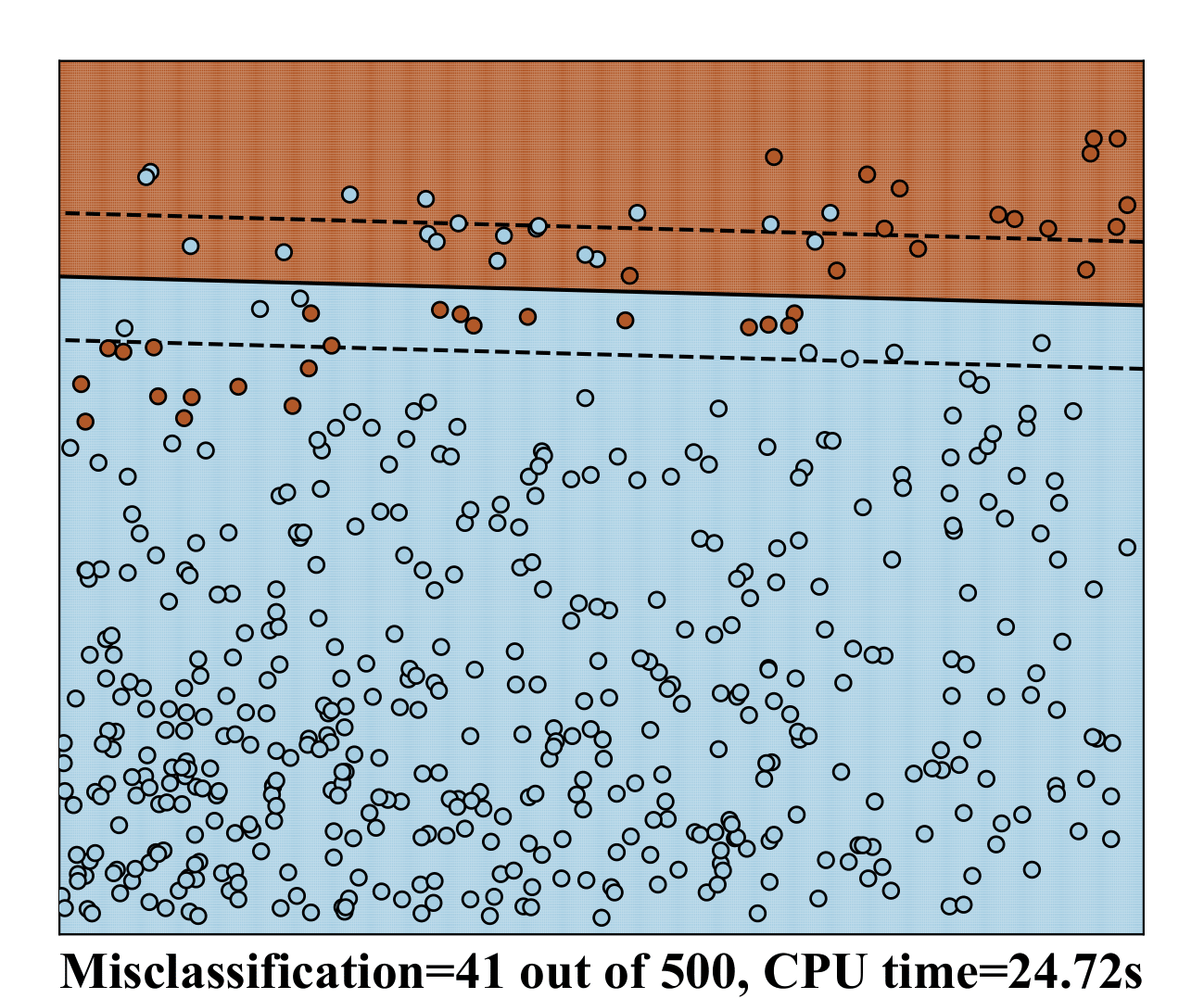} &
        \includegraphics[width=4.1cm]{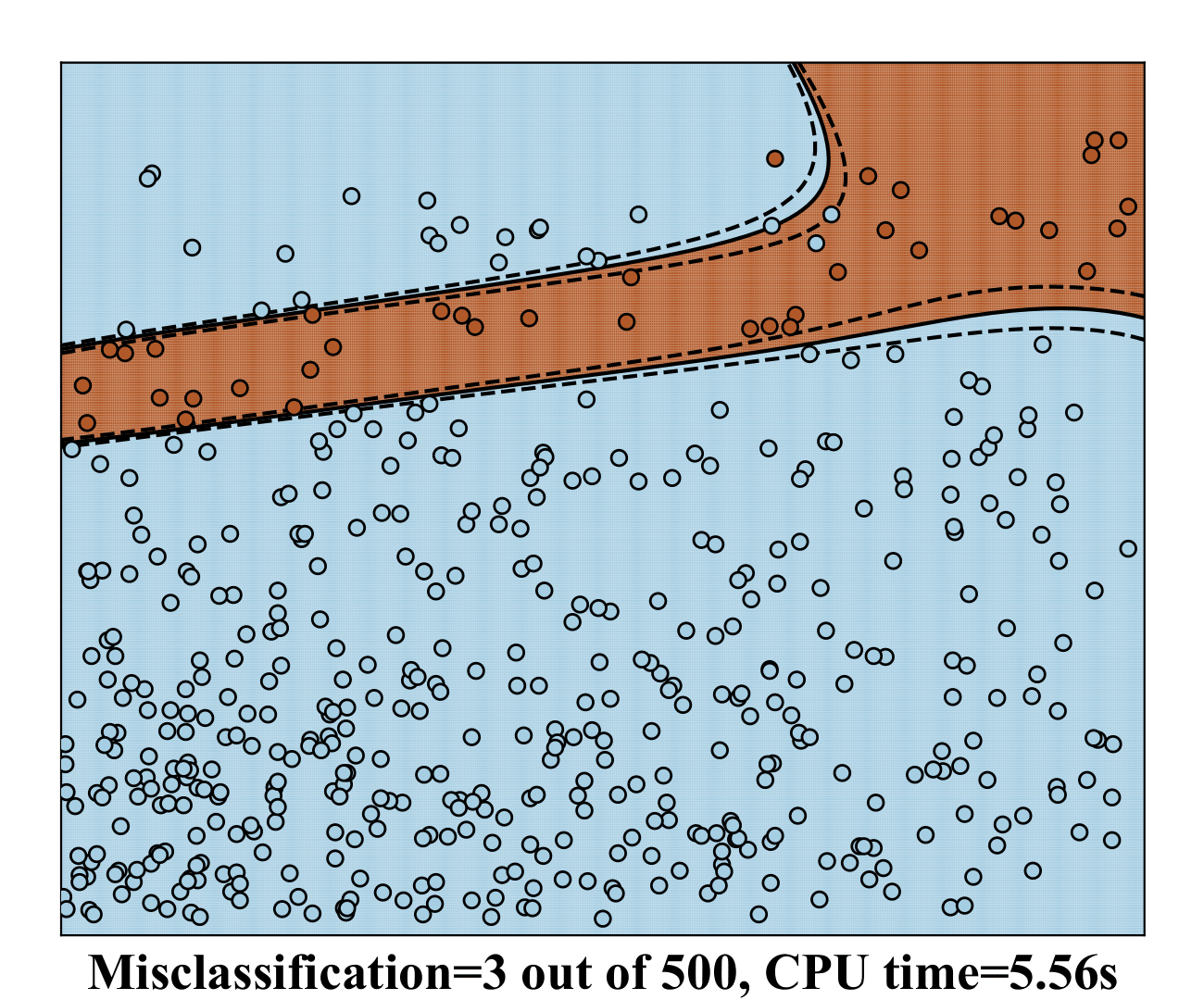} &
        \includegraphics[width=4.1cm]{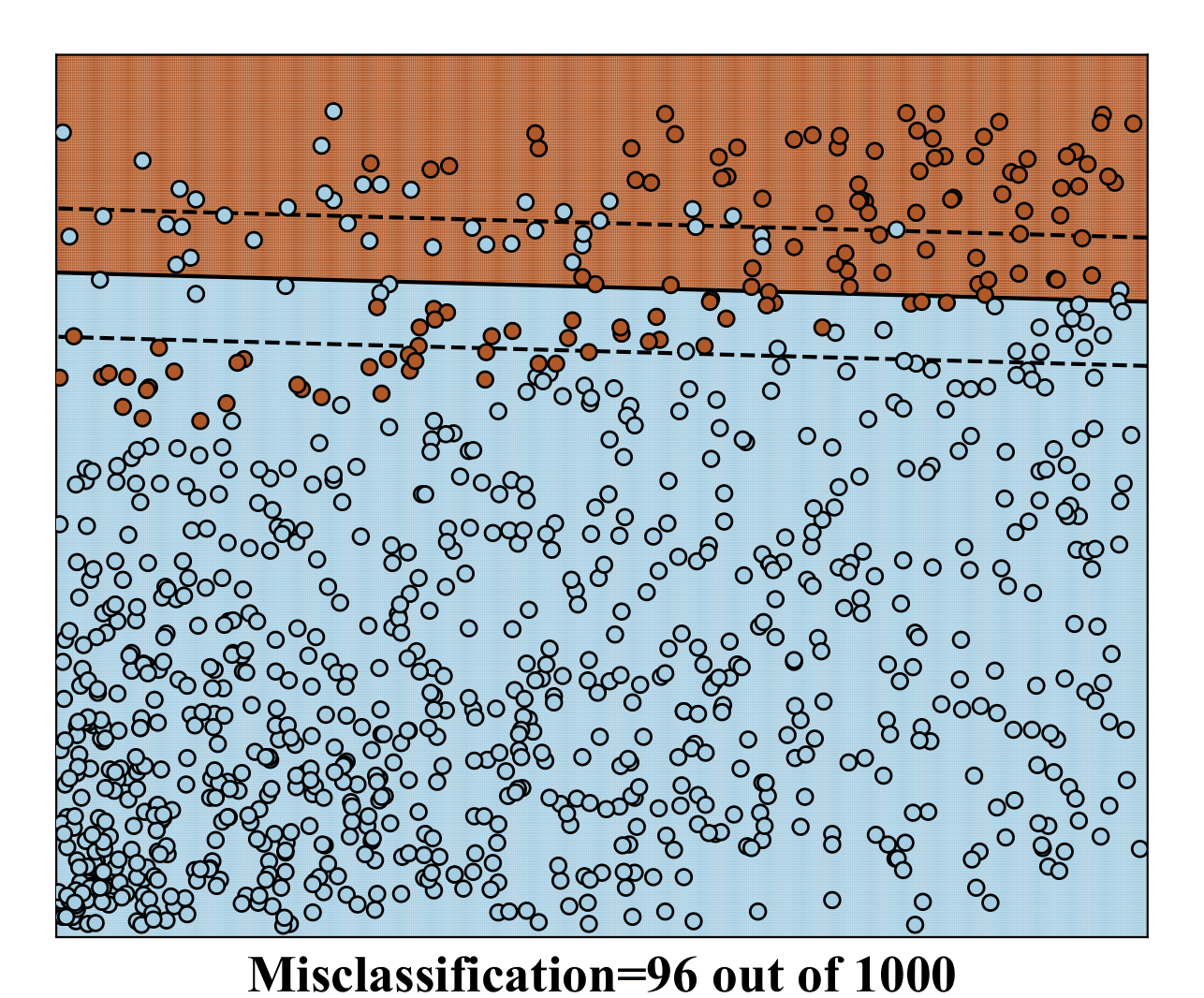} &
        \includegraphics[width=4.1cm]{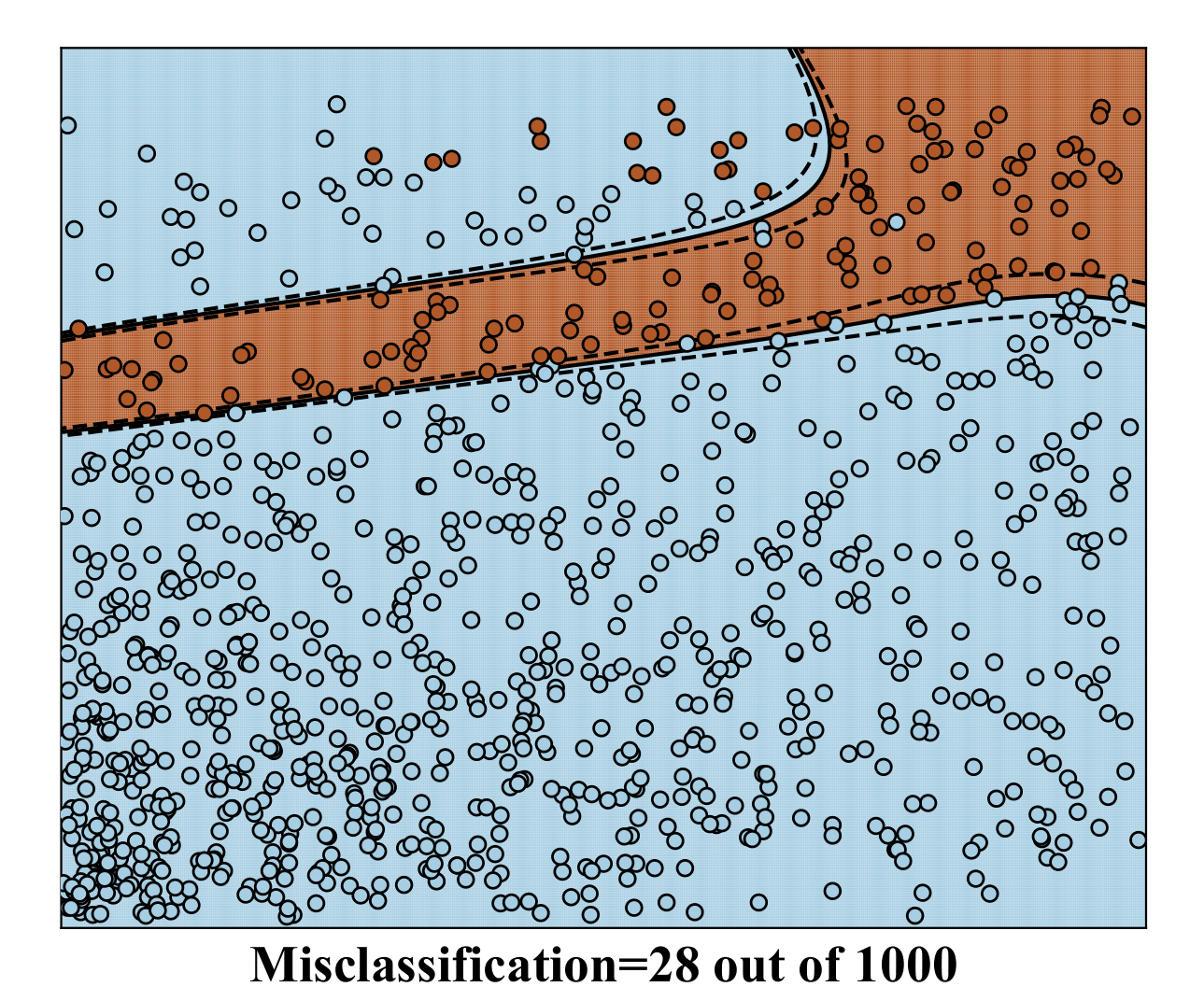} \\
        \includegraphics[width=4.1cm]{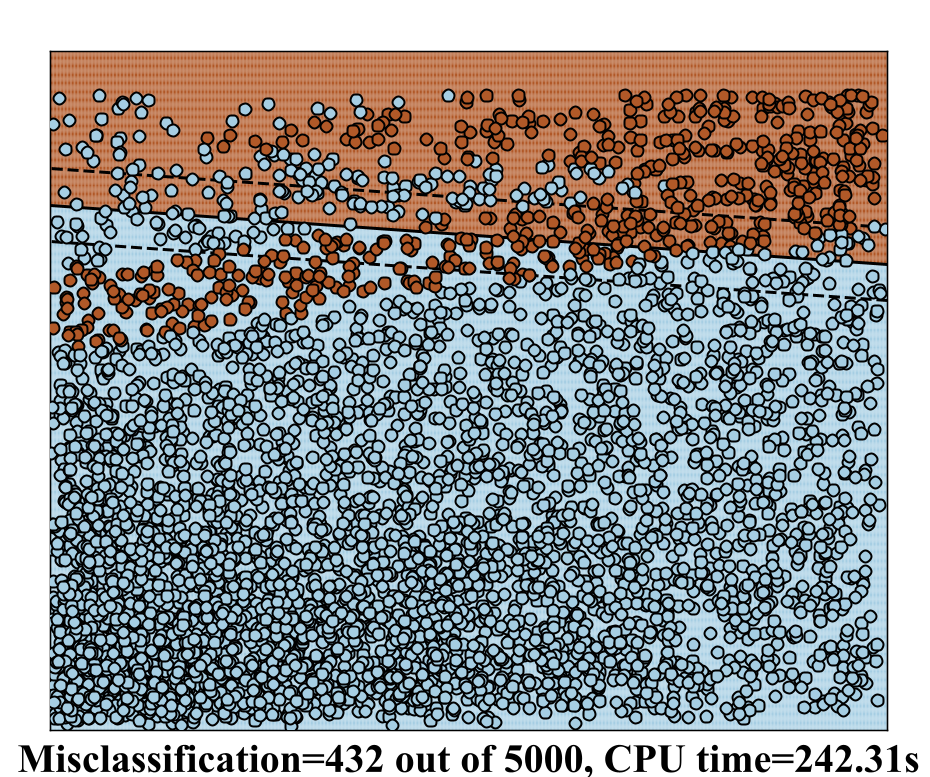} &
        \includegraphics[width=4.1cm]{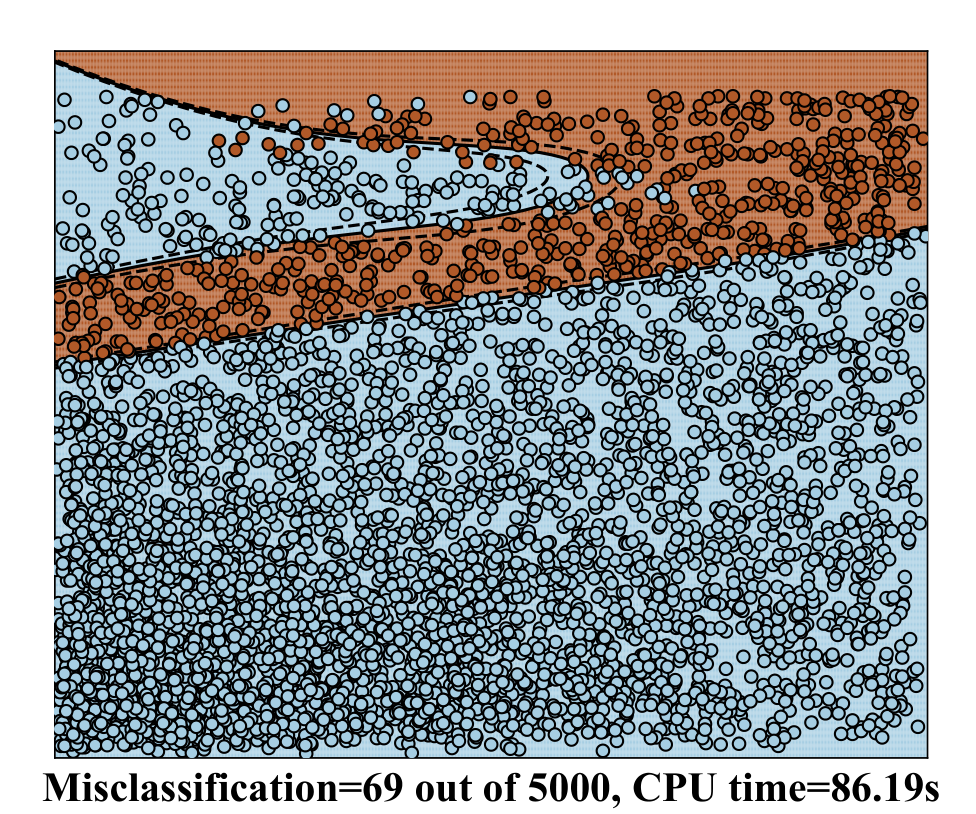} &
        \includegraphics[width=4.1cm]{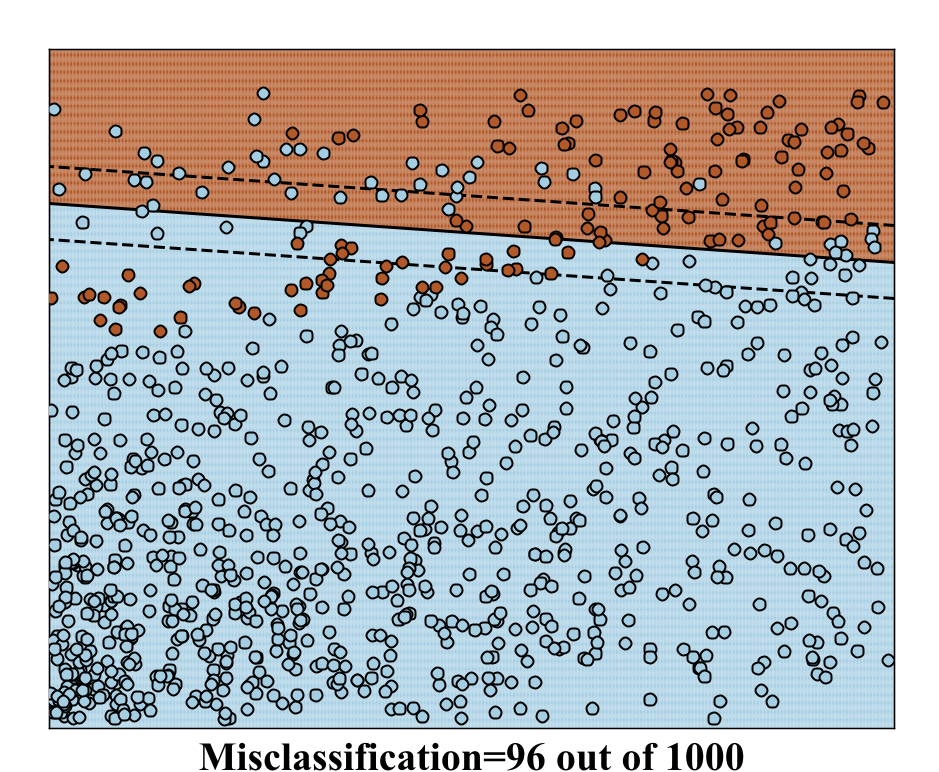} &
        \includegraphics[width=4.1cm]{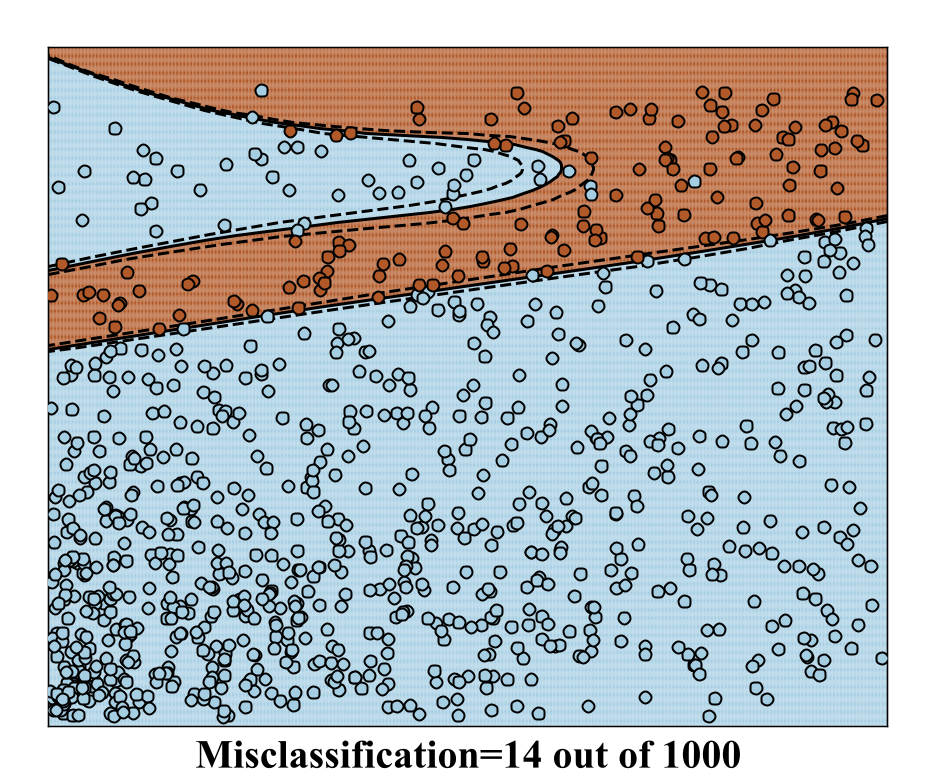} \\\\
       };
\node [white, font=\bfseries] at (-8.0,4.3) {$\textbf{(1)}$};
\node [white, font=\bfseries] at (-3.7,4.3) {$\textbf{(2)}$};
\node [white, font=\bfseries] at (+0.6,4.3) {$\textbf{(3)}$};
\node [white, font=\bfseries] at (+4.9,4.3) {$\textbf{(4)}$};
\node [white, font=\bfseries] at (-8.0,0.7) {$\textbf{(5)}$};
\node [white, font=\bfseries] at (-3.7,0.7) {$\textbf{(6)}$};
\node [white, font=\bfseries] at (+0.6,0.7) {$\textbf{(7)}$};
\node [white, font=\bfseries] at (+4.9,0.7) {$\textbf{(8)}$};
\node [white, font=\bfseries] at (-8.0,-3.0) {$\textbf{(9)}$};
\node [white, font=\bfseries] at (-3.6,-3.0) {$\textbf{(10)}$};
\node [white, font=\bfseries] at (+0.7,-3.0) {$\textbf{(11)}$};
\node [white, font=\bfseries] at (+5.0,-3.0) {$\textbf{(12)}$};
\end{tikzpicture}
\caption{The number of missclassified samples and CPU time for SVM classifiers using $50$, $500$, and $5,000$ training samples and $1,000$ test samples.}
\label{result1}
\vspace{-2mm}
\end{figure*}

Once $\tilde{\mathbf{y}}$ predicted, we check the feasibility and the upper bound guarantee. If the prediction is feasible and the guarantee is desirable, the prediction is directly adopted as the commitment decision for conventional units. Otherwise, we use $\tilde{\mathbf{y}}$ as a warm start to the MISOCP solver built upon a branch-and-bound algorithm. Note that for this, we convert back $\tilde{y}_{gt}=-1$ to $\tilde{y}_{gt}=0$, indicating the corresponding unit is off. 
 The branch-and-bound algorithm uses the warm start as an initial solution. By that, it uses the warm start as a bound to eliminate parts of the search space that cannot produce a better solution and  therefore focus on the rest of search space. 
This is particularly useful when the MISOCP problem has a large number of variables, and the solution space is vast. We find a trade-off between prediction quality and computational time to determine how much mixed-integer solvers can benefit from utilizing the predicted $\tilde{\mathbf{y}}$ as a warm start.

\color{black}
Recall that the predicted $\tilde{\mathbf{y}}$ to be used as a warm start might be infeasible for the MISOCP problem. 
If $\tilde{\mathbf{y}}$ is infeasible since it fails to meet the minimum up- and down-time constraints, we find a close feasible solution among the sampled strategies to the infeasible one by using a heuristic approach to recover a feasible solution. The heuristic approach sorts 
 the training dataset $(\mathbf{X,Y})$ based on increasing distances $||\mathbf{y}_h-\Tilde{\mathbf{y}}||_2$ for all $h=\{1,...,H\}$. Then, it selects the nearest $\mathbf{y}_h$ to $\tilde{\mathbf{y}}$ represented by $\hat{\mathbf{y}}$ to be used as a warm start. In the case where there are multiple $\hat{\mathbf{y}}$ having the same distance to $\tilde{\mathbf{y}}$, we can select the one that occurs the most in the training dataset. Knowing $\mathbf{y}_h$ for all $h=\{1,...,H\}$ satisfies the minimum up- and down-time constraints, the warm start $\hat{\mathbf{y}}$ meets these constraints as well. Furthermore, in order to guarantee that there is always a feasible solution to the UC problem in terms of satisfying power flow constraints, we also take extreme corrective measures such as load shedding and renewable power curtailment into consideration. 

\color{black}

\section{Experimental results}
We provide two case studies based on the IEEE $6$-bus and $118$-bus test systems. We solve the optimization problems in Python using scikit-learn API and in Matlab using the YALMIP toolbox with Gurobi solver $8.1.1$ on a $16$ GB RAM personal computer clocking at $3.1$ GHZ. \color{black} All source codes are publicly available\footnote{\url{https://github.com/farzanehpourahmadi/UCP.git}}. \color{black}

\subsection{IEEE 6-bus test system}
This stylized system includes $3$ conventional units, $2$ wind farms, and $3$ loads, as shown in Fig. \ref{6bus}. The technical data for conventional units and transmission lines is given in \cite{GitHub}. For simplicity in this stylized case, 
we assume that loads are deterministic and that two wind farms only have stochastic generation. Additionally, we assume an identical shape for daily wind power generation profile across different samples --- the only varying factor is their generation level. As such, we end up in two features only, i.e., the production level of two wind farms. Wind and load data is provided in \cite{GitHub}.

\subsubsection*{Training data} We utilize three distinct sets of training data for our analysis, which contain $50$, $500$, and $5,000$ samples, respectively. Recall for generating each sample $(\mathbf{x},\mathbf{y})$, we solve the MISOCP problem given in the appendix with the feature vector $\mathbf{x}$ as the input data, and compute commitments $\mathbf{y}$. Based on Algorithm 1 in Section III, given $\epsilon=0.10$, these number of training samples provide a desired probability guarantee of $\delta=1$, $\delta=0.40$, and $\delta=0.18$, respectively. This means that, if one provides further samples beyond $50$, $500$, and $5,000$ training samples, the maximum probability of observing a new strategy with the confidence level of $90\%$ is $100\%$, $40\%$, and $18\%$, respectively. Having samples $(\mathbf{x},\mathbf{y})$ as inputs data, we solve the linear or conic optimization problems for SVM classifiers (see Table \ref{classfiers}), and obtain map functions \eqref{map1} for the linear SVM and \eqref{map2} for the kernelized SVM.

\subsubsection*{Testing data} We use a dataset containing $1,000$ samples in the form of $(\tilde{\mathbf{x}},\tilde{\mathbf{y}})$, which are different than training samples. We use the map functions  \eqref{map1} and \eqref{map2} to predict $\tilde{\mathbf{y}}$, to be used either as a commitment decision or to serve for warm starting. 

\subsubsection*{Cross-validation for kernel function and regularizer} We have used a cross-validation approach to tune the parameter $\gamma$ in the Gaussian kernel function. \color{black} To do so, we split training samples into $4$ equally sized subsets and then training the model on $3$ subsets, while validating its performance on the remaining subset. This process is repeated $4$ times, with each subset being used as the validation set once, and the other $3$ subsets being used for training. The performance metric used is the expected hinge loss of validation instances. Similarly, we regularize both linear and kernelized SVM classifiers by determining a value for $\lambda$, obtained by cross-validation too. Note that
we employ cross-validation on the training dataset, without using the test dataset, to prevent a biased training process.

In order to justify the selection of four folds, we have conducted a sensitivity analysis regarding the number of folds used in cross-validation. We have varied the number of folds from $3$ to $6$ and observed that there is no notable change in the expected hinge loss. Specifically, we have found that the change is less than $0.1\%$ for both linear and kernelized SVM classifiers with regularization. By this, we have concluded $4$-fold is a reasonable choice for cross-validation. \color{black}  

\subsubsection*{Classification results} Let us first explore how successful the map functions \eqref{map1} and \eqref{map2} are in the prediction of the correct commitments. Fig. \ref{result1} shows classification boundaries for the commitment of  arbitrarily selected conventional unit $\rm{G}3$ in a certain hour. For every plot, the two axes are wind production of two farms (two available features). Based on classification results, it is predicted  $\rm{G}3$ to be off in the blue area and to be on in the brown area. The blue and brown circles indicate the \textit{true} commitment, obtained by solving the MISOCP problem in the appendix. If a blue (brown) circle is  located in the blue (brown) area, it shows that the prediction is correct, otherwise it is a misclassification.
This figure includes $12$ plots, where those in the first, second, and third rows correspond to cases where the number of training samples is $50$, $500$, and $5,000$, respectively. The plots in the first two columns pertain to the training stage, where those plots in the third and fourth columns show the testing results with  $1,000$ samples. Under each plot, the number of misclassified samples and the corresponding computational time (in the training stage) are reported.


One main observation is that the kernelized SVM outperforms the linear one in terms of the number of misclassified samples in both training and testing stages. For example, in the case with $500$ training samples, the kernelized SVM in the training stage (plot $6$) has only $3$ out of $500$ samples wrongly classified, whereas it is $41$ for the linear SVM (plot $5$). Similarly, in the testing stage, the kernelized SVM misclassified $28$ out of $1,000$ samples (plot $8$), while it is $96$ samples for the linear SVM (plot $7$). There is a similar observation in the case with $5,000$ training samples, i.e., plots $9$-$12$. This has been expected as the kernelized SVM offers a higher degree of freedom when determining a map function, resulting in nonlinear boundaries between blue and brown areas. The interesting point is that the computational time for training the kernelized SVM is lower than the linear one, making it  even a more appealing choice. Another interesting observation is that by increasing the number of training samples from $500$ to $5,000$, the kernelized SVM exhibits a more successful performance in the testing stage (reduced misclassified samples from $28$ in plot $8$ to $14$ in plot $12$). However, it is not the case for the linear SVM (see plots $7$ and $11$).

\subsubsection*{Bias analysis and performance guarantee}
\begin{figure}
    \color{black}\centering
    \includegraphics[width=1\columnwidth]{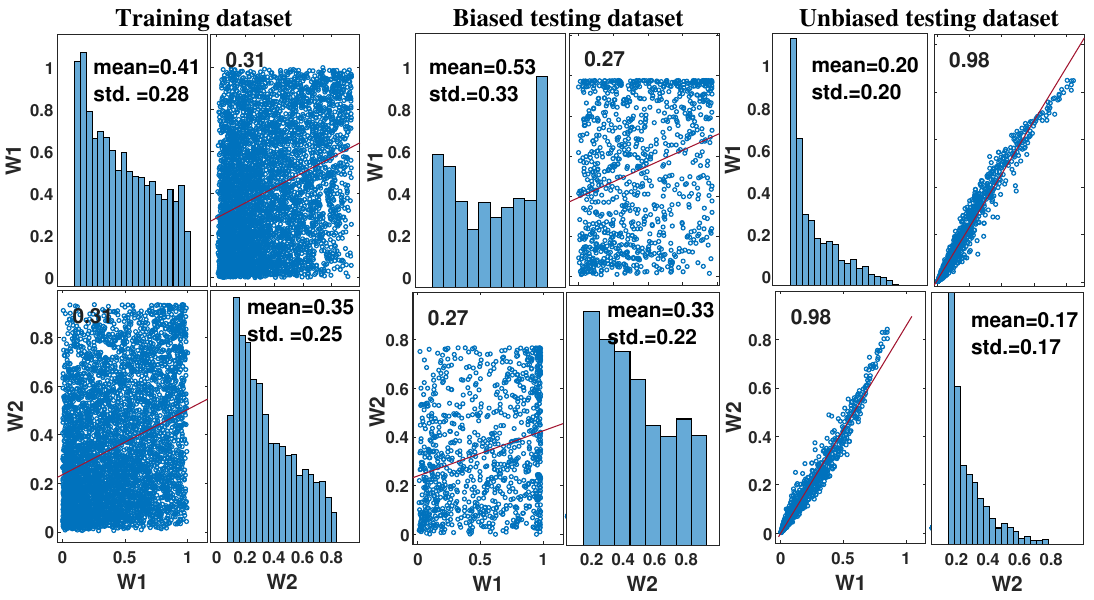}
    \caption{A matrix of plots for the training dataset and biased and unbiased testing datasets. The diagonal plots (x-axis: in per-unit; y-axis: frequency of occurrence) display the wind production level frequency.
The off-diagonal plots (both axes in per-unit) display the correlation of two wind farms $\rm{W}1$-$\rm{W}2$.}
    \label{correlation}
    \vspace{-1mm}
\end{figure}

We first start with a bias analysis in terms of the similarity of training and test datasets. We keep the previous training dataset with $5,000$ samples. However, we generate two different test datasets, each with $1,000$ samples, one is more similar to the training dataset, while the other one is less similar. We use the  $2$-Wasserstein metric to measure the distance of training and test datasets. The similar test dataset, so called \textit{biased} test dataset, has a distance of $0.05$ to the training samples. In contrast, the other one, so called \textit{unbiased} test dataset, has a comparatively higher distance of $0.20$. \color{black} Fig. \ref{correlation} presents a matrix of plots for the training and two biased and unbiased testing datasets separately. Diagonal plots illustrate wind production level histograms, whereas scatter plots depict the correlation between the two wind farms. As observed in Fig. \ref{correlation}, the biased testing dataset closely resembles the unbiased one in terms of mean, standard deviation (std.), and correlation coefficients of the two wind farms. \color{black} We then exploit the kernelized SVM classifier \eqref{dual_kernel_primal} including a regularization term. 

\begin{figure}
    \centering
    \includegraphics[width=0.9\columnwidth]{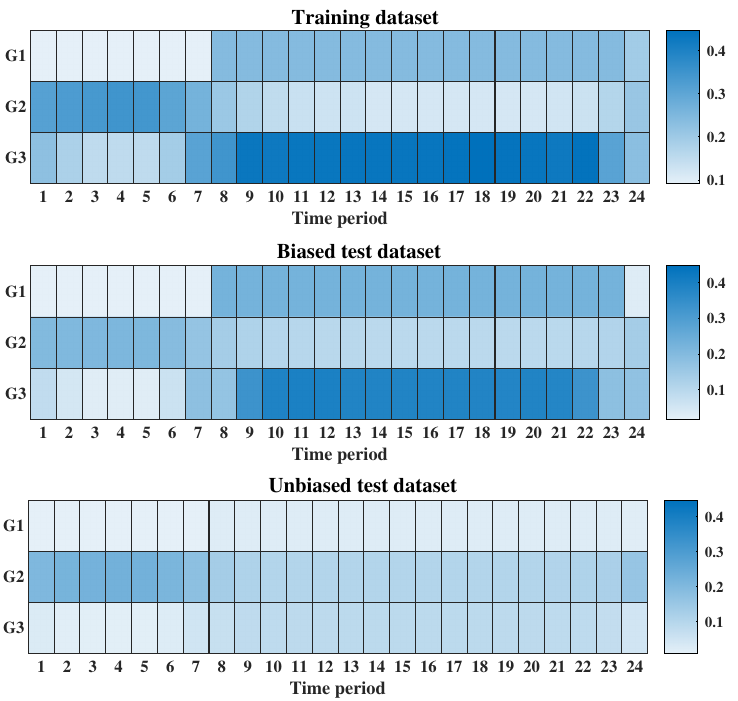}
    \caption{The expected hinge loss of the kernelized SVM classifier, for every unit and hour, in the training and testing stages (with biased and unbiased test data).}
    \label{heatmap_loss}
    \vspace{-1mm}
\end{figure}

For all three conventional units over $24$ hours, Fig. \ref{heatmap_loss} reports the optimal value of the expected hinge loss of the kernelized SVM classifier in both training and testing stages. For the training dataset, recall that the expected hinge loss is equal to \eqref{dual_kernel_primal} in the optimal point, whereas for the test datasets, it is computed as the left-hand side of \eqref{predictionguarantee_kernel}. Owed to the regularization weight $\lambda$ carefully tuned by cross validation, the expected hinge loss in the unbiased test dataset is lower than that in the biased one, which is desirable. In addition, the expected hinge loss in both biased and unbiased test datasets is lower than that of the training dataset, validating the performance guarantee \eqref{predictionguarantee_kernel} is obtained. 

\color{black}
\subsubsection*{Impact of initial states}
\begin{figure}
\color{black}
    \centering
    \includegraphics[width=1\columnwidth]{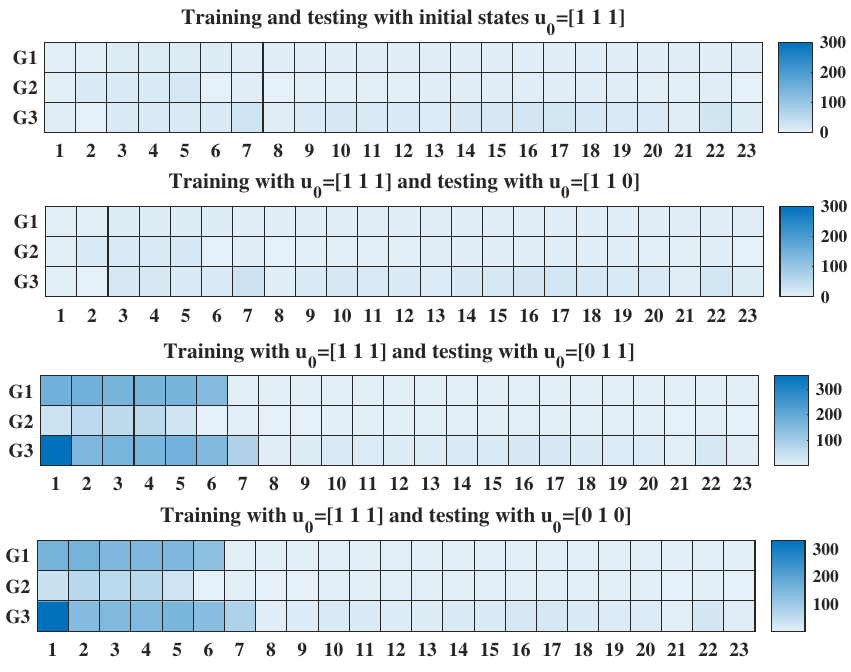}
    \caption{The accuracy of predictions when the initial states $\mathbf{u}_0$ for testing is different. Accuracy here refers to the percentage of incorrect predictions out of $1,000$ testing samples.}
    \label{initialstates5}
    \vspace{-1mm}
\end{figure}

We now investigate the impact of initial states of conventional  units on the accuracy of their commitments' predictions. So far, we have not restricted the initial on/off states of $\rm{G}1$-$\rm{G}3$ in our test system. However, by introducing constraints on the initial states of units in a UC problem, we can compile a dataset that includes the impact of these initial states. Let us suppose that initially,  units $\rm{G}1$-$\rm{G}3$ are all on, represented by $\mathbf{u}_0=[1,1,1]$. We derive the corresponding UC solution $\mathbf{y}$ accordingly to construct our training dataset. Subsequently, we train binary classifiers using this dataset. Fig. \ref{initialstates5} shows how accurately the prediction reflects the initial states $\mathbf{u}_0$ considered in the testing dataset, indicating whether they are similar or different from the initial states $\mathbf{u}_0$ considered in the training dataset.

The upper plot in Fig. \ref{initialstates5} illustrates the accuracy of predictions when the initial states $\mathbf{u}_0$ for testing match with those used in constructing the training dataset. Accuracy here refers to the percentage of incorrect predictions out of $1,000$ testing samples. Additionally, we examine three other scenarios where the initial states $\mathbf{u}_0$ of testing samples differ from those of the training set. We consider four cases only in Fig. \ref{initialstates5} as  unit $\rm{G}2$ remains on owed to its lowest operating cost among three units. Therefore, there are $4$ possible combinations of initial states for the two remaining units.

From Fig. \ref{initialstates5}, we observe that if the initial state of $\rm{G}3$ changes to $0$ during testing, the prediction accuracy remains unaffected. However, for $\rm{G}1$, it is crucial to train binary classifiers specifically for the corresponding initial states. This observation is also valid for larger test systems. Initially, one needs to identify the set of potential initial states and determine which ones require retraining of binary classifiers. It is important to note that all these investigations can be conducted offline. 

\color{black}
Similarly, for any change in the optimization formulation of the UC problem, one may
need retraining the classifier. This could be the case when a conventional unit or a transmission line is out for maintenance, or when the system operator changes the network topology.

\color{black}
\subsubsection*{Benchmark model} We have found models developed in \cite{Juanmi} and \cite{xavier} relevant and comparable to the UC predictor proposed in this work. They employ a classical supervised learning method called k-nearest neighbors (KNN) to predict the commitment of conventional units. We adopt KNN as a benchmark for constructing binary classifiers and compare its prediction performance with that of the proposed SVM classifiers. Furthermore, \cite{Juanmi} concludes that the performance of existing and upcoming learning-based methods to solve the UC problem should be compared with some simple and straightforward, yet effective methods like KNN. Inspired by this viewpoint, we treat KNN as a benchmark and investigate the effectiveness of both KNN and SVM classifiers in predicting the correct commitments.

The rationale behind KNN is to identify the most similar examples in the training dataset and conduct a majority vote. Specifically, for a given feature vector $\Tilde{\mathbf{x}}$, the method sorts the training dataset $(\mathbf{X,Y})$ based on increasing distances $||\mathbf{x}_h-\Tilde{\mathbf{x}}||_2$ for all $h=\{1,...,H\}$. It selects the first $k$ nearest neighbors to $\Tilde{\mathbf{x}}$ with the lowest distances, denoted as $\{(\mathbf{x}_1,\mathbf{y}_1),...,(\mathbf{x}_k,\mathbf{y}_k)\}$. Then, to predict each binary variable $\Tilde{y}_{gt}$, it computes the average of the values $\frac{y_{1gt}+...+y_{kgt}}{k}$. If the average is positive, $\Tilde{y}_{gt}=1$; otherwise, $\Tilde{y}_{gt}=-1$. It is worth noting that similar to SVM, we map $\Tilde{y}_{gt}=-1$ back to $\Tilde{y}_{gt}=0$.

\begin{figure}
    \color{black}\centering
    \includegraphics[width=0.95\columnwidth]{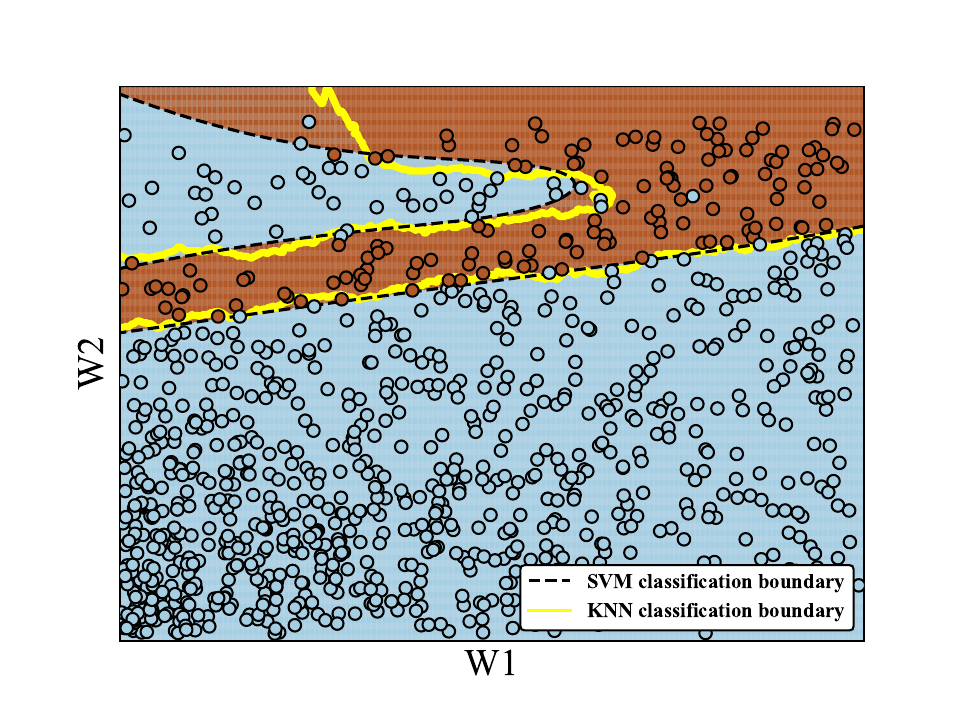}
    \caption{Comparison of Kernelized SVM with regularization and KNN classifiers.}
    \label{KNN_SVM}
    \vspace{-1mm}
\end{figure}

Consider plot $12$ in Fig. \ref{result1}, illustrating the classification boundary generated by Kernelized SVM for the commitment of a randomly chosen conventional unit $\rm{G}3$ during a specific hour. We then apply the KNN technique to this particular example to determine its classification boundary. It is worth mentioning that we have employed a cross-validation approach to fine-tune the parameter $k$ in the KNN technique. A comparison between SVM and KNN boundaries is depicted in Fig. \ref{KNN_SVM}. In both cases, $14$ out of $1,000$ test samples are misclassified. We repeat this process for other binary classifiers corresponding to different  units and hours, observing that SVM and KNN exhibit almost identical performance. 
This observation confirms the conclusion in \cite{Juanmi}, stating that KNN might work satisfactory for the UC classification such that more complex classifiers might be unnecessary. However, we realize that it holds true for this small-scale $6$-bus test system. In the next section with a larger case study, we will observe that the proposed SVM classifier outperforms KNN by a system cost saving of  $5\%$ when there is a time limit for the computation task or  a computational time reduction of $3\%$ when there is no such a limit.
\color{black}

\subsection{IEEE 118-bus test system}
We apply the proposed UC predictor to the IEEE $118$-bus test system, including  $19$ conventional units, $2$ wind farms, $91$ loads, and $186$ transmission lines. Our main focus is on the potential improvement in the computational performance. Furthermore, we investigate the potential cost saving by improving the optimality gap\footnote{\black{We realize the computational requirements are way more significant when solving MISOCP problem for cases larger than the $118$-bus test system. In the relevant UC literature including AC power flow equations \cite{UC_AC1,UC_AC2,UC_AC3,UC_AC4,UC_AC5}, we found out that the largest system considered in those studies is of a comparable scale to the $118$-test system.}}.
\color{black}

\subsubsection*{Data preparation}
We use a dataset comprising of $5,000$ samples for training and $1,000$ samples for testing. 
The features that we consider here include the power generation of wind farms as well as the load in all buses over $24$ hours. Therefore, each classifier for unit $g$ and hour $t$ has access to $93$ features. Note that we assume all loads and wind farms adhere to a fixed power factor.  Additionally, we assume an identical daily load and wind power generation profile across various samples. Times series for daily load and wind power generation are given in \cite{GitHub}. For each sample, we solve the MISOCP problem (see appendix), including $2,592$ binary variables, $13,920$ continuous variables, and $55,464$ linear and conic constraints. With a confidence level of $1-\epsilon=0.90$, Algorithm 1 in Section III validates a  probability guarantee of $\delta=0.12$ to observe a new strategy. This means, with the confidence level of $90$\%, there is a maximum probability of $12$\% to observe a new strategy if further training samples are provided. 
Similar to the previous case study, we tune the parameters for regularization and Gaussian kernel function using a cross-validation approach with $4$ equally sized subsets.

\begin{figure*}[t] \label{out-of-sample}
\centering
\begin{minipage}{0.5\textwidth}
  \centering
  	\scalebox{1}{
	\includegraphics[width=1\columnwidth]{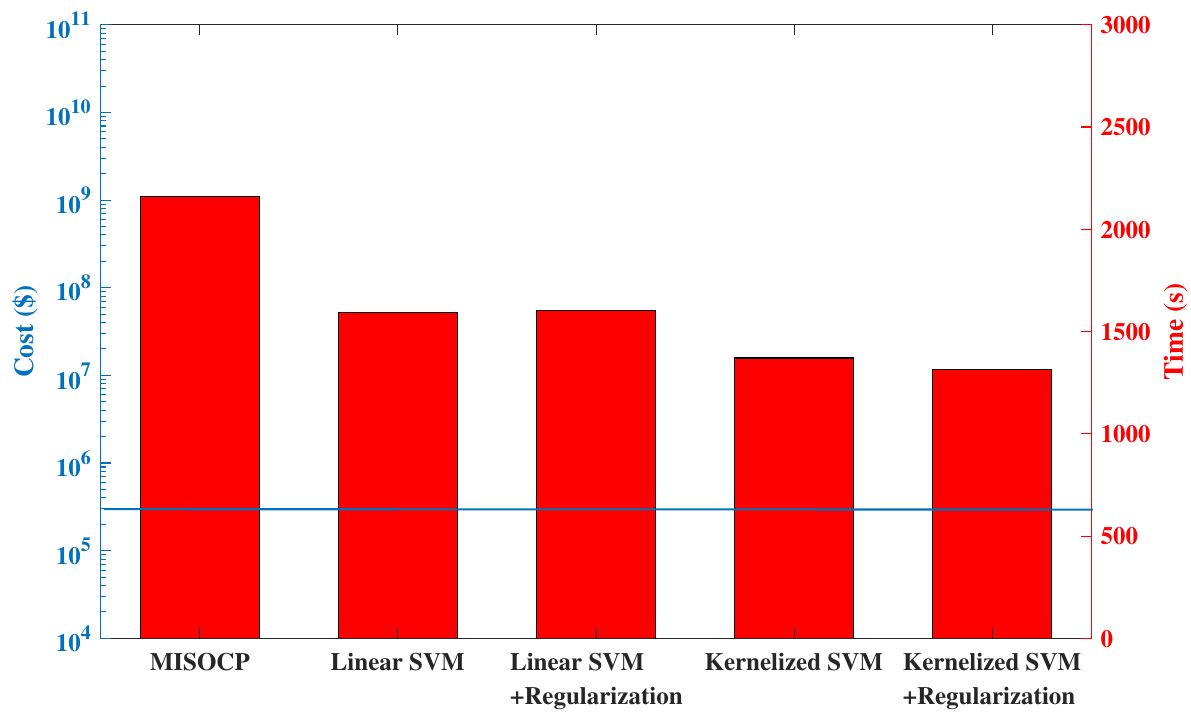}
	}
\end{minipage}%
\begin{minipage}{.5\textwidth}
  \centering
    	\scalebox{0.98}{
	\includegraphics[width=1\columnwidth]{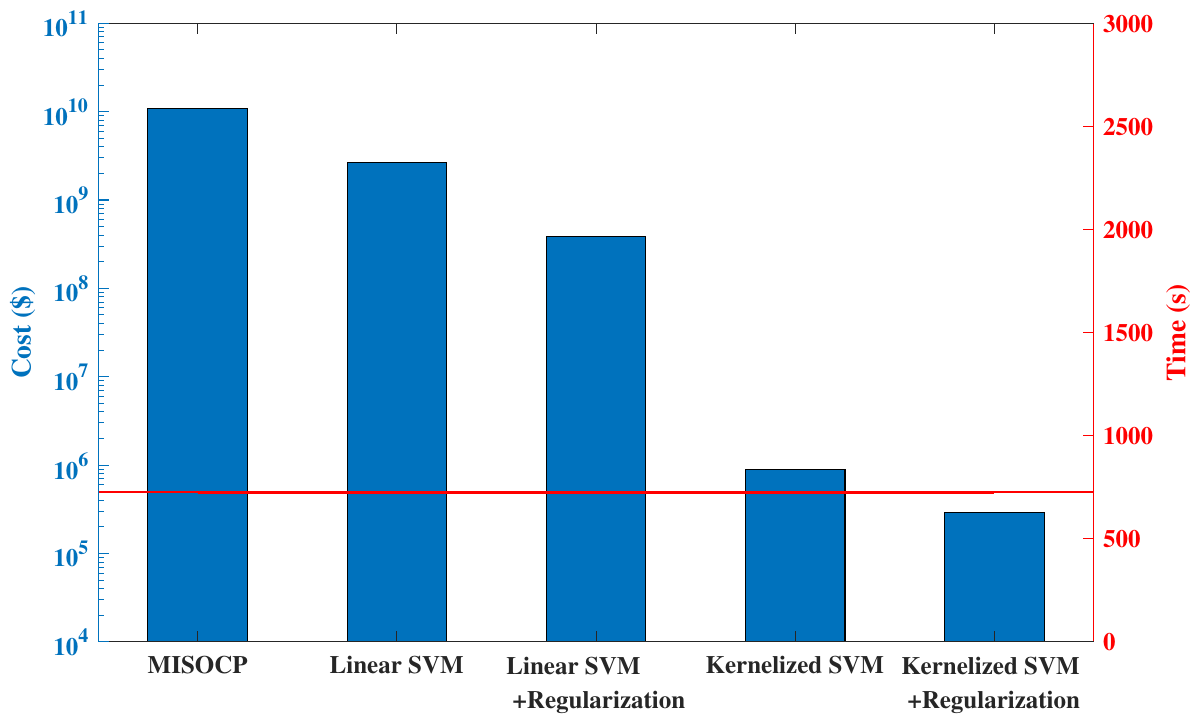}
	}
\end{minipage}
    \caption{Out-of-sample performance: Average system cost and average computational time for solving the UC problem, solved by the MISOCP solver  without (bar 1) or with warm start (bars 2-5). If warmly started, it is initiated by the prediction of classifiers. The right plot enforces a computational time limit of $700$ seconds, and shows the system cost obtained in such a limited period. The left plot shows the computational time when no limit is enforced.}
\label{time limitation}
    \vspace{-2mm}
\end{figure*}

\subsubsection*{Implementation}
We train a classifier for every unit and hour. 
If the predicted strategy $\tilde{\mathbf{y}}$ is desirable, we can directly use it to determine the commitment of conventional units, but we must still verify that the minimum up- and down-time requirements are satisfied.
Otherwise, if the predicted strategy $\tilde{\mathbf{y}}$ is not desirable, we leverage it as a warm start for the Gurobi solver and explore whether it speeds up the computational time to solve the UC problem, i.e., the MISOCP problem in the appendix. \color{black} If $\tilde{\mathbf{y}}$ is not binary feasible, we find a feasible strategy that is as close as possible to the predicted strategy among the strategies sampled for training dataset. \color{black} The distance between strategies can be measured as $\lVert \widehat{\mathbf{y}} - \tilde{\mathbf{y}}\rVert$, where $\tilde{\mathbf{y}}$ is the predicted strategy and $\widehat{\mathbf{y}}$ is the closest feasible strategy, respectively.

\subsubsection*{Evaluation}

\begin{figure}
    \centering
    \includegraphics[width=0.95\columnwidth]{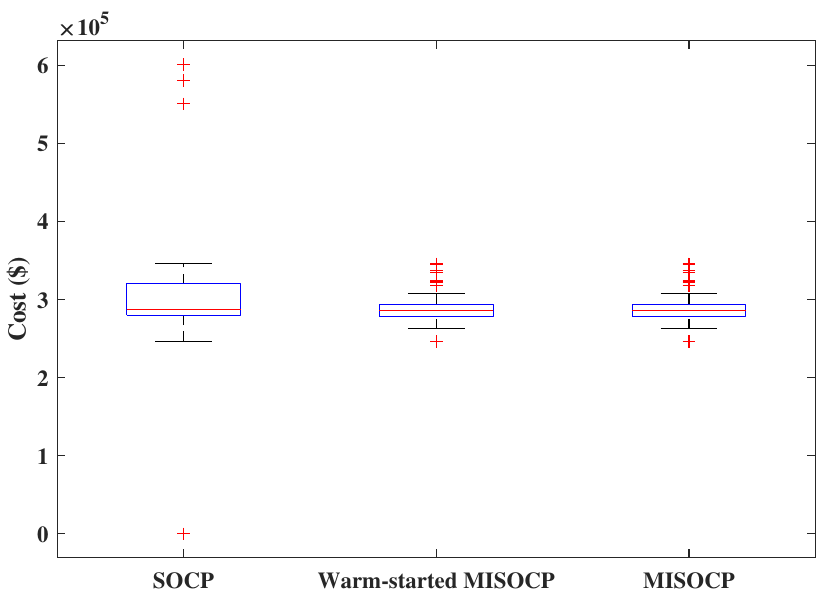}
    \caption{The median and the $25$th and $75$th percentiles of the out-of-sample system cost obtained by solving the UC problem. First bar: UC problem with fixed predictions (binary variables) resulting in a second-order cone programming (SOCP) problem. Second and third bars: UC problem as a MISOCP problem with and without warm start.}
    \label{SOCP}
    \vspace{-2mm}
\end{figure}


In the setting that we use the prediction as a warm start to the MISOCP problem, it is likely to end up in a sub-optimal prediction, but it may still help reduce the computational time. 
Hereafter, we acknowledge a UC prediction to be \textit{effective}  if it enhances the computational performance of the UC solver. Therefore, our measurement for the prediction effectiveness is to what extent it reduces the computational time. In the following, we examine two scenarios: in Scenario I, the system operator has a limit for the computational time which cannot be exceeded. This implies that the UC problem may not be solved to optimality if the computational time reaches to its limit. In Scenario II, there is not such a limit, and therefore the UC problem can be solved to optimality. 

Now, we answer two questions: 
(\textit{i}) What is the least possible difference between the optimal strategies $\mathbf{y}^{*}$ and the predicted ones $\tilde{\mathbf{y}}$ under Scenario I (with time limit)?
(\textit{ii}) How quickly can we obtain the optimal strategies $\mathbf{y}^{*}$ under Scenario II (no time limit) with warm starting? To answer these two questions, we train four different classifiers, namely the linear SVM (with and without regularization) and the kernelized SVM (with and without regularization). We then compare the UC results with and without using warm starts --- in the case of warm start, the classifiers provide the initial solution.

In the following, with or without warm start, we solve the UC problem $1,000$ times (each time with a different test dataset), and report the average out-of-sample system cost and the average computational time. By the out-of-sample system cost, we refer to the total generation and start-up cost of the system, as formulated in the objective function \eqref{objective} of the UC problem in the appendix.

The right plot of Fig. \ref{time limitation} answers the first question under scenario I, where the computational time limit of $700$ seconds has been enforced. This implies when the computational time has reached to this limit, we stopped the solver, and reported the system cost.  One can expect that this cost is significantly high if the solver has not managed to find a solution with a low optimality gap in the limited computational period.  The right plot of Fig. \ref{time limitation} shows that the warmly started MISOCP problem initiated by the predictions of the kernelized SVM with regularization (last bar) has been solved effectively, resulting in the lowest system cost on average. \color{black} On the contrary, the MISOCP solver without warm start (first bar) is very slow, such that it could not solve the UC problem effectively in $700$ seconds, yielding a high optimality gap, and thereby a significant system cost. We observe that a significant portion of the running time is usually spent in finding a feasible solution. Therefore, the MISOCP solver warmly started by the predictions of other classifiers exhibits a performance in between. \color{black} In general, this plot concludes that the kernelized SVM provides more effective warm starts than the linear one. 

\color{black}
In cases where there is no computational time constraint or a warm start has been provided, we have not observed any uncertainty in the output of the solver --- by no uncertainty, we mean the solver provides the same solution in every run. However, in the case a computational time limit exists without any warm start, as illustrated in the first bar on the left plot of Fig. \ref{time limitation}, we have observed an uncertainty. Nevertheless, our general conclusion remains valid that utilizing SVM-enabled warm starts reduces the system cost, despite the existence of such an uncertainty in certain situations.
\color{black}

The left plot of Fig. \ref{time limitation} answers the second question, under scenario II where there is no computational time limit. The MISOCP solver has obtained the same UC solution irrespective of warm or cold start. However, the computational time is significantly higher without warm start. The kernelized SVM with regularization is again the best choice for predicting an effective warm start, reducing the computational time by an average factor of $1.7$, compared to a case with no warm start. 

We now evaluate the out-of-sample system cost if we trust in predicted strategies $\tilde{\mathbf{y}}$ and directly use them as commitment decisions of conventional units. We fix the binary variables in the UC problem to those predictions, resulting in a SOCP problem. Fig. \ref{SOCP} shows the median and the $25$th and $75$th percentiles of the out-of-sample system cost when we solve the SOCP problem with fixed binary variables, and compares it to those results of MISOCP without or with warm start (by the kernelized SVM with regularization). It has been observed that the system cost when fixing binary variables is $5$\% higher than that of the MISOCP problem, but it is $4.5$ times faster in terms of computational time than the MISOCP problem without warm start, and $2.4$ times faster than the same problem with warm start.

\color{black}
\subsubsection*{Benchmark model}
We employ KNN predictions as a warm start for the MISOCP problem. Similar to our approach in evaluating different SVM techniques, we explore two scenarios to assess KNN's prediction effectiveness in reducing computational costs. In scenario I, constrained by a computational time limit of $700$ seconds, the system cost using KNN predictions as a warm start is $\$277,010$. 
The kernelized SVM continues to provide more effective warm starts, resulting in a $5\%$ improvement compared to KNN. In scenario II, where there are no computational time constraints, utilizing KNN predictions as a warm start allows solving the UC problem to optimality with an expected computational time of $1,354$ seconds. Warm starting with kernelized SVMs reduces the computational time by $3\%$ compared to KNN.
\color{black}

\section{Conclusion}
This paper showed that, if the system operators have restricted by computational time limits to solve the UC problem, there is a high potential to use SVM-based classification methods to predict commitment decisions of conventional units. As a pragmatic solution, these predictions can be used as a warm start to considerably speed up the mixed-integer optimization solvers, without compromising the optimality. We showed that it is possible to properly regularize the SVM classifiers in a way that it provides an out-of-sample performance guarantee. In addition, this paper discussed how many training samples are needed to get an insight into the probability of having unseen strategies. This paper concluded that the kernelized SVM with a proper regularization (tuned by cross validation) outperforms the linear SVM, and can help the system operators to ease their computational burden. 

This paper has developed one classifier per conventional unit per hour. This overlooks the potential spatio-temporal correlations, which could be very useful information when training the classifiers. It is of interest to develop such correlation-aware classifiers, exploring to what extent it may provide more effective predictions (e.g., for warm start) and helps to further reduce the computational time. Another future direction is to explore the performance of classifiers in non-stationary environments, where one may not learn much for old data. 
\color{black} Another direction for future work is to explore how to incorporate the convex hull description of UC decisions proposed in \cite{Kai}, aiming to enhance the computational efficiency of the UC problem as an expansion of this study. Additionally, an intriguing direction for further investigation is to explore methods for transfer learning to adapt to changes in the UC problem, such as modifications in network topology, without a need for retraining the binary classifier.
\color{black}

\section*{Appendix: MISOCP Formulation for UC Problem}
%
\textit{Notation}: We consider a power grid consisting of $N$ buses, $L$ transmission lines, and $M$ conventional units. For every hour $t$, the set of continuous variables contains active power generation of conventional units $\mathbf{p}^{\rm{G}}_t=\{p^{\rm{G}}_{1t},...,p^{\rm{G}}_{gt},...,p^{\rm{G}}_{Mt}\}\in \R^M$,  reactive power generation of conventional units  $\mathbf{q}^{\rm{G}}_t=\{q^{\rm{G}}_{1t},...,q^{\rm{G}}_{gt},...,q^{\rm{G}}_{Mt}\}\in \R^M$, nodal voltage magnitudes $\boldsymbol{\nu}_t=\{\nu_{1t},...,\nu_{it},...,\nu_{Nt}\}\in \R^N$, nodal voltage angles $\boldsymbol{\theta}_t=\{\theta_{1t},...,\theta_{it},...,\theta_{Nt}\}\in \R^N$, and apparent power flow of transmission lines $\mathbf{s}_{t}=\{s_{1t},...,s_{lt},...,s_{Lt}\}\in \R^L$. In addition, the set of binary variables includes on/off commitment of conventional units $\mathbf{u}_t=\{u_{1t},...,u_{gt},...,u_{Mt}\}\in \{0,1\}^M$ and their start-up status $\mathbf{v}_t=\{v_{1t},...,v_{gt},...,v_{Mt}\}\in \{0,1\}^M$. 

We now define parameters. For each conventional unit $g$, parameters $\overline p_g/\underline p_g$ and $\overline{q}_g/\underline q_g$ indicate the maximum/minimum active and reactive power generation limits, respectively. Parameters $\overline{r}_g/\underline r_g$ denote the maximum ramp-up and ramp-down capabilities, whereas $\tilde{r}_g$ gives the ramp-rate capability in start-up and shut-down hours. Parameter $\underline{v}_g$ gives the minimum up- and down-time, respectively. In addition, $\mathbf{c}^P=\{c^P_{1},...,c^P_{g},...,c^P_{M}\}\in \R^M$ and $\mathbf{c}^S=\{c^S_{1},...,c^S_{g},...,c^S_{M}\}\in \R^M$ give the production and
start-up costs of conventional units, respectively. For each transmission line $l$ connecting bus $i$ to bus $j$, $G_{ij}$ and $B_{ij}$ denote the real and imaginary parts of the bus admittance matrix, respectively. In addition, $B^{\rm{sh}}_{ij}$ represents the reactive shunt element, considering a $\pi$-model for transmission lines. The apparent power capacity of line $l$ is given by $\overline{s}_l$. For each bus $i$, $\overline{\nu}_i/\underline{\nu}_i$ gives the maximum/minimum voltage magnitude. Parameters $\mathbf{p}^{\rm{W}}_t=\{p^{\rm{W}}_{1t},...,p^{\rm{W}}_{it},...,p^{\rm{W}}_{Nt}\}\in \R^N$ and $\mathbf{q}^{\rm{W}}_t=\{q^{\rm{W}}_{1t},...,q^{\rm{W}}_{it},...,q^{\rm{W}}_{Nt}\}\in \R^N$ denote active and reactive power generation of wind farms, respectively. Finally, $\mathbf{p}^{\rm{D}}_t=\{p^{\rm{D}}_{1t},...,p^{\rm{D}}_{it},...,p^{\rm{D}}_{Nt}\}\in \R^N$ and $\mathbf{q}^{\rm{D}}_t=\{q^{\rm{D}}_{1t},...,q^{\rm{D}}_{it},...,q^{\rm{D}}_{Nt}\}\in \R^N$ represent nodal active and reactive power consumptions, respectively.  

The original AC unit commitment problem reads as 
\begingroup
\allowdisplaybreaks
\begin{subequations}
\begin{align}
&\min_{\mathbf{p},\mathbf{q},\mathbf{u},\mathbf{v},\boldsymbol{\nu},\boldsymbol{\theta}} \ \ \sum_t  \mathbf{p}^{\rm{G}}_t (\mathbf{c}^{P})^\top+ \mathbf{v}_t (\mathbf{c}^{S})^\top\label{objective} \\
    &-\!u_{g(t-1)}\!+\!u_{gt}-u_{g\tau}\!\leq\! 0, \forall \tau\! \in\! \{t,...,\!\overline{v}_g+t\!-\!1\}, \forall g,t \label{minup}\\
    &\;u_{g(t-1)}-u_{gt}+u_{g\tau}\leq\! 1,  \forall \tau\! \in\! \{t,...,\!\underline{v}_g+t\!-\!1\}, \forall g,t \label{mindn}\\
    &-u_{g(t-1)}+u_{gt}-v_{gt}\leq 0, \enspace \forall g,t \label{statetran} \\
    &u_{gt}\underline{p}_g\leq p^{\rm{G}}_{gt} \leq \overline{p}_g u_{gt}, \enspace \forall g,t \label{active}\\
    &u_{gt}\underline{q}_g\leq q^{\rm{G}}_{gt} \leq \overline{q}_g u_{gt}, \enspace \forall g,t \label{reactive} \\
    &p^{\rm{G}}_{gt}-p^{\rm{G}}_{g(t-1)}\leq \overline{r}_g u_{g(t-1)}+\tilde{r}_g(1-u_{g(t-1)}), \enspace \forall g,t \label{ramp1}\\
    &p^{\rm{G}}_{g(t-1)}-p^{\rm{G}}_{gt}\leq \underline{r}_g u_{gt}+\tilde{r}_g(1-u_{gt}), \enspace \forall g,t \label{ramp2}\\
    &\sum_{g\in \Gamma(i)}\!\!p^{\rm{G}}_{gt}+p^{\rm{W}}_{it}-p^{\rm{D}}_{it}=G_{ii}\nu_{it}^2+ 
    \!\sum_{j\in \Lambda(i)}\!\!\!\nu_{it}\nu_{jt} \nonumber\\
    &\hspace{15mm}\Big(G_{ij}\cos(\theta_{it}\!-\!\theta_{jt})\!\!+\!B_{ij}\sin(\theta_{it}\!-\!\theta_{jt})\Big), \forall i,t \! \label{balance_active}\\
    &\sum_{g\in \Gamma(i)}\!\!q^{\rm{G}}_{gt}+q^{\rm{W}}_{it}-q^{\rm{D}}_{it}=-B_{ii}\nu_{it}^2+ \!\sum_{j\in \Lambda(i)}\!\!\!\nu_{it}\nu_{jt} \nonumber\\
    &\hspace{15mm}\Big(G_{ij}\sin(\theta_{it}\!-\!\theta_{jt})\!\!-\!B_{ij}\cos(\theta_{it}\!-\!\theta_{jt})\Big), \forall i,t \label{balance_reactive}\\
    & s^2_{lt}=\Big(\!-G_{ij}\nu_{it}^2\!+\!G_{ij}\nu_{it}\nu_{jt}\cos(\theta_{it}\!-\!\theta_{jt})\!\nonumber\\
    &\hspace{5mm}+\!B_{ij}\nu_{it}\nu_{jt}\sin(\theta_{it}\!-\!\theta_{jt})\Big)^2+\Big(\!(B_{ij}-B_{ij}^{\rm{sh}})\nu_{it}^2\!\nonumber\\
    &\hspace{5mm}-\!B_{ij}\nu_{it}\nu_{jt}\cos(\theta_{it}\!-\!\theta_{jt})\!+\!G_{ij}\nu_{it}\nu_{jt}\sin(\theta_{it}\!-\!\theta_{jt})\Big)^2, \nonumber\\
    &\hspace{60mm}\forall l\in (i,j),t\label{Equation}\\
    & s^2_{lt}\leq \overline{s}_l^2, \enspace  \forall l\in (i,j),t \label{line}\\
   & \underline{\nu}_i\le \nu_{it}\leq \overline{\nu}_i, \enspace \forall i,t \label{voltage} \\
   & \theta_{ref,t}=0, \enspace\forall t \label{angleref}\\
   &\textbf{u}_t,\textbf{v}_t \in \{0,1\} ,\enspace\forall t. \label{lastconstraint}
\end{align}
\end{subequations}
\endgroup

The objective function (\ref{objective}) minimizes the total production and start-up costs of conventional units. Constraints (\ref{minup})-(\ref{mindn}) restrict the minimum up- and down-time, whereas (\ref{statetran}) determines the state transition of conventional units. Capacity and ramping limits of conventional units are enforced by (\ref{active})-(\ref{ramp2}). Constraints (\ref{balance_active})-(\ref{balance_reactive}) represent the nodal active and reactive power balance, where $\Gamma(i)$ and $\Lambda(i)$ denote the set of conventional units connected to bus $i$ and the set of adjacent buses to  $i$, respectively. Constraint (\ref{Equation}) determines the apparent power flow in line $(i,j)$, constrained by (\ref{line}). Constraint (\ref{voltage}) limits the voltage magnitudes. Constraint (\ref{angleref}) sets voltage angle of the reference bus to zero. Finally, (\ref{lastconstraint}) declares the binary variables. 

This AC unit commitment includes non-convex constraints (\ref{balance_active})-(\ref{Equation}), as they are quadratic equality constraints. For the convexification of those constraints, the current literature \cite{jabr} suggests defining new variables for each bus $i$ as $w_{it}=\nu_{it}^2$ and for each line $(i,j)$ as $e_{ijt}=\nu_{it}\nu_{jt}\cos(\theta_{it}-\theta_{jt})$ and $f_{ijt}=-\nu_{it}\nu_{jt}\sin(\theta_{it}-\theta_{jt})$. This lets replace (\ref{balance_active})-(\ref{line}) by
\begingroup
\allowdisplaybreaks
\begin{subequations}
\begin{align}
&\sum_{g\in \Gamma(i)}\!\!\!p^{\rm{G}}_{gt}\!+\!p^{\rm{W}}_{it}\!\!-\!p^{\rm{D}}_{it}\!= \nonumber\\
&\hspace{2cm}\!G_{ii}\nu_{it}\!+\!\!\!\!\sum_{j\in \Lambda(i)}\!\!\!(G_{ij}e_{ijt}\!\!-\!\!B_{ij}f_{ijt}), \forall i,\!t \label{balance_active_re}\\
&\sum_{g\in \Gamma(i)}\!\!\!q^{\rm{G}}_{gt}\!+\!q^{\rm{W}}_{it}\!\!-\!q^{\rm{D}}_{it}\!= \nonumber\\
&\hspace{2cm}\!-B_{ii}\nu_{it}\!-\!\!\!\!\!\!\sum_{j\in \Lambda(i)}\!\!\!(B_{ij}e_{ijt}\!\!+\!\!G_{ij}f_{ijt}), \forall i,\!t\!  \label{balance_reactive_re}\\
& \Big(\!-G_{ij}\nu_{it}\!+\!G_{ij}e_{ijt}-\!B_{ij}f_{ijt}\Big)^2+\!\!\Big(\!(B_{ij}\!-\!B_{ij}^{\rm{sh}})\nu_{it}\!\nonumber\\
&\hspace{16mm}-\!B_{ij}e_{ijt}\!-\!G_{ij}f_{ijt}\Big)^2\!\!\!\le\! \overline{s}_l^2, \forall l\!\in\! (i,j),\!t\label{line_ref}\\
&0=e_{ijt}^2+f_{ijt}^2-\nu_{it}\nu_{jt},\forall i,j,t \label{efu}\\
& 0=\theta_{jt}-\theta_{it}-{\rm{atan}}(\frac{f_{ijt}}{e_{ijt}}),\forall i,j,t \label{angel_ref}\\
&\underline{v}_i^2\le \nu_{it}\le \overline{v}_i^2,\forall i,t \label{u_ref}\\
&-\overline{v}_i\overline{v}_j \le e_{ijt},f_{ijt} \le \overline{v}_i\overline{v}_j,\forall i,j,t. \label{ef_ref}
\end{align}
\end{subequations}
\endgroup

Now, the original formulation converts to (\ref{objective})-(\ref{ramp2}), (\ref{voltage})-(\ref{lastconstraint}),  (\ref{balance_active_re})-(\ref{ef_ref}), which is still non-convex. However, it can be relaxed to a mixed-integer second-order cone optimization by omitting (\ref{angel_ref}) and relaxing (\ref{efu}) as 
\begin{align}
  e_{ijt}^2+f_{ijt}^2\le \nu_{it}\nu_{jt},\forall i,j,t. \label{efu_relax}  
\end{align}

The resulting MISOCP problem is (\ref{objective})-(\ref{ramp2}), (\ref{voltage})-(\ref{lastconstraint}),  (\ref{balance_active_re})-(\ref{line_ref}), (\ref{u_ref})-(\ref{ef_ref}) and (\ref{efu_relax}), which is used to solve the UC problem in this paper.

\section*{Acknowledgment}
We would like to thank Thomas Falconer (DTU) for reading the manuscript and providing constructive feedback. We also thank three anonymous reviewers for constructive feedback. 




\bibliographystyle{IEEEtran}
\bibliography{IEEEabrv,ML}

\vfill
	
	\markboth{ }{ }
	\setcounter{figure}{0}
	\setcounter{table}{0}
	
	
	\onecolumn
	
	\begingroup
	\begin{large}
	%
	\setlength{\parskip}{1em}
	\setlength{\parindent}{0em}

	\end{large}
	\endgroup

\end{document}

%% file: comments2.tex
\usepackage{blindtext}
\usepackage{multirow}
\usepackage{bbm}
\usepackage{dsfont}
\usepackage[scr]{rsfso}
\usepackage{cite}
\usepackage[cmex10]{amsmath}
\usepackage{amssymb}
\usepackage{bbm}
\usepackage{dsfont}
\usepackage{mathtools}
\usepackage{color}
\usepackage{amsthm}
\usepackage{mathabx}
\usepackage{lipsum}

\newtheorem{theorem}{Theorem}

\newcommand{\be}{\begin{equation}}
\newcommand{\ee}{\end{equation}}
\newcommand{\bea}{\begin{equation*}\begin{aligned}}
\newcommand{\eea}{\end{aligned}\end{equation*}}

\newcommand{\R}{\mathbb{R}}




%% file: main.bbl
\begin{thebibliography}{10}
\providecommand{\url}[1]{#1}
\csname url@samestyle\endcsname
\providecommand{\newblock}{\relax}
\providecommand{\bibinfo}[2]{#2}
\providecommand{\BIBentrySTDinterwordspacing}{\spaceskip=0pt\relax}
\providecommand{\BIBentryALTinterwordstretchfactor}{4}
\providecommand{\BIBentryALTinterwordspacing}{\spaceskip=\fontdimen2\font plus
\BIBentryALTinterwordstretchfactor\fontdimen3\font minus
  \fontdimen4\font\relax}
\providecommand{\BIBforeignlanguage}[2]{{%
\expandafter\ifx\csname l@#1\endcsname\relax
\typeout{** WARNING: IEEEtran.bst: No hyphenation pattern has been}%
\typeout{** loaded for the language `#1'. Using the pattern for}%
\typeout{** the default language instead.}%
\else
\language=\csname l@#1\endcsname
\fi
#2}}
\providecommand{\BIBdecl}{\relax}
\BIBdecl

\bibitem{anjos}
M.~F. Anjos and A.~J. Conejo, ``Unit commitment in electric energy systems,''
  \emph{Foundations and Trends® in Electric Energy Systems}, vol.~1, no.~4,
  pp. 220--310, 2017.

\bibitem{hobbs}
B.~F. Hobbs, M.~H. Rothkopf, R.~P. O'Neill, and H.-p. Chao, \emph{The next
  generation of electric power unit commitment models}.\hskip 1em plus 0.5em
  minus 0.4em\relax Springer Science \& Business Media, 2006, vol.~36.

\bibitem{UC_AC1}
J.~Liu, C.~D. Laird, J.~K. Scott, J.-P. Watson, and A.~Castillo, ``Global
  solution strategies for the network-constrained unit commitment problem with
  ac transmission constraints,'' \emph{IEEE Trans. Power Syst.}, vol.~34,
  no.~2, pp. 1139--1150, 2019.

\bibitem{tejada}
D.~A. Tejada-Arango, P.~S{\'a}nchez-Mart{\i}n, and A.~Ramos, ``Security
  constrained unit commitment using line outage distribution factors,''
  \emph{IEEE Trans. Power Syst.}, vol.~33, no.~1, pp. 329--337, 2017.

\bibitem{linear}
C.~Coffrin and P.~Van~Hentenryck, ``A linear-programming approximation of {AC}
  power flows,'' \emph{INFORMS J. Comput.}, vol.~26, no.~4, pp. 718--734, 2014.

\bibitem{SOCP}
R.~Madani, M.~Ashraphijuo, and J.~Lavaei, ``Promises of conic relaxation for
  contingency-constrained optimal power flow problem,'' \emph{IEEE Trans. Power
  Syst.}, vol.~31, no.~2, pp. 1297--1307, 2015.

\bibitem{SDP}
C.~Coffrin, H.~L. Hijazi, and P.~Van~Hentenryck, ``Strengthening the {SDP}
  relaxation of {AC} power flows with convex envelopes, bound tightening, and
  valid inequalities,'' \emph{IEEE Trans. Power Syst.}, vol.~32, no.~5, pp.
  3549--3558, 2016.

\bibitem{atakan}
S.~Atakan, G.~Lulli, and S.~Sen, ``A state transition {MIP} formulation for the
  unit commitment problem,'' \emph{IEEE Trans. Power Syst.}, vol.~33, no.~1,
  pp. 736--748, 2017.

\bibitem{morales}
G.~Morales-Espa{\~n}a, J.~M. Latorre, and A.~Ramos, ``Tight and compact {MILP}
  formulation for the thermal unit commitment problem,'' \emph{IEEE Trans.
  Power Syst.}, vol.~28, no.~4, pp. 4897--4908, 2013.

\bibitem{rajan}
D.~Rajan and S.~Takriti, ``Minimum up/down polytopes of the unit commitment
  problem with start-up costs,'' \emph{IBM Res. Rep}, vol. 23628, pp. 1--14,
  2005.

\bibitem{ostrowski}
J.~Ostrowski, M.~F. Anjos, and A.~Vannelli, ``Tight mixed integer linear
  programming formulations for the unit commitment problem,'' \emph{IEEE Trans.
  Power Syst.}, vol.~27, no.~1, pp. 39--46, 2011.

\bibitem{tightening1}
D.~A. Tejada-Arango, S.~Lumbreras, P.~Sánchez-Martín, and A.~Ramos, ``Which
  unit-commitment formulation is best? {A} comparison framework,'' \emph{IEEE
  Trans. Power Syst.}, vol.~35, no.~4, pp. 2926--2936, 2020.

\bibitem{tightening2}
C.~Gentile, G.~Morales-España, and A.~Ramos, ``A tight {MIP} formulation of
  the unit commitment problem with start-up and shut-down constraints,''
  \emph{EURO J. Comput. Optim.}, vol.~5, no.~1, pp. 177--201, 2017.

\bibitem{Kai}
K.~Pan and Y.~Guan, ``Convex hulls for the unit commitment polytope,''
  \emph{arXiv preprint arXiv:1701.08943}, 2017.

\bibitem{decomposition}
Y.~Fu, M.~Shahidehpour, and Z.~Li, ``Security-constrained unit commitment with
  ac constraints,'' \emph{IEEE Trans. Power Syst.}, vol.~20, no.~2, pp.
  1001--1013, 2005.

\bibitem{costscreening}
Q.~Zhai, X.~Guan, J.~Cheng, and H.~Wu, ``Fast identification of inactive
  security constraints in scuc problems,'' \emph{IEEE Trans. Power Syst.},
  vol.~25, no.~4, pp. 1946--1954, 2010.

\bibitem{Cheol}
D.~Bertsimas and C.~W. Kim, ``A prescriptive machine learning approach to
  mixed-integer convex optimization,'' \emph{INFORMS J. Comput.}, vol.~35,
  no.~6, pp. 1225--1241, 2023.

\bibitem{D1}
F.~Fioretto, T.~W. Mak, and P.~Van~Hentenryck, ``Predicting {AC} optimal power
  flows: Combining deep learning and {L}agrangian dual methods,'' in
  \emph{Proceedings of the AAAI Conference on Artificial Intelligence},
  vol.~34, 2020, pp. 630--637.

\bibitem{Pascal_PSCC}
W.~Chen, S.~Park, M.~Tanneau, and P.~Van~Hentenryck, ``Learning optimization
  proxies for large-scale security-constrained economic dispatch,''
  \emph{Electr. Power Syst. Res.}, vol. 213, p. 108566, 2022.

\bibitem{Proxy_UC1}
G.~Dalal, E.~Gilboa, S.~Mannor, and L.~Wehenkel, ``Unit commitment using
  nearest neighbor as a short-term proxy,'' in \emph{2018 Power Systems
  Computation Conference (PSCC)}, 2018, pp. 1--7.

\bibitem{Proxy_UC2}
X.~Lin, Z.~J. Hou, H.~Ren, and F.~Pan, ``Approximate mixed-integer programming
  solution with machine learning technique and linear programming relaxation,''
  in \emph{2019 3rd International Conference on Smart Grid and Smart Cities
  (ICSGSC)}, 2019, pp. 101--107.

\bibitem{ML1}
D.~Bertsimas and B.~Stellato, ``Online mixed-integer optimization in
  milliseconds,'' \emph{INFORMS J. Comput.}, vol.~34, no.~4, p. 2229–2248,
  2022.

\bibitem{roaldimplied}
L.~Roald and D.~Molzahn, ``Implied constraint satisfaction in power system
  optimization: The impacts of load variations,'' in \emph{Annual Allerton
  Conference on Communication, Control, and Computing (Allerton)}.\hskip 1em
  plus 0.5em minus 0.4em\relax IEEE, 2019, pp. 308--315.

\bibitem{ID1}
D.~Deka and S.~Misra, ``Learning for {DC-OPF: C}lassifying active sets using
  neural nets,'' in \emph{IEEE Milan PowerTech}, 2019, pp. 1--6.

\bibitem{ID2}
S.~Misra, L.~Roald, and Y.~Ng, ``Learning for constrained optimization:
  Identifying optimal active constraint sets,'' \emph{INFORMS J. Comput.},
  vol.~34, no.~1, pp. 463--480, 2022.

\bibitem{Juanmi}
S.~Pineda and J.~Morales, ``Is learning for the unit commitment problem a
  low-hanging fruit?'' \emph{Electr. Power Syst. Res.}, vol. 207, p. 107851,
  2022.

\bibitem{UC_ActiveConstraint1}
S.~Pineda, J.~M. Morales, and A.~Jiménez-Cordero, ``Data-driven screening of
  network constraints for unit commitment,'' \emph{IEEE Trans. Power Syst.},
  vol.~35, no.~5, pp. 3695--3705, 2020.

\bibitem{UC_ActiveConstraint2}
D.~T. Lagos and N.~D. Hatziargyriou, ``Data-driven frequency dynamic unit
  commitment for island systems with high {RES} penetration,'' \emph{IEEE
  Trans. Power Syst.}, vol.~36, no.~5, pp. 4699--4711, 2021.

\bibitem{xavier}
{\'A}.~S. Xavier, F.~Qiu, and S.~Ahmed, ``Learning to solve large-scale
  security-constrained unit commitment problems,'' \emph{INFORMS J. Comput.},
  vol.~33, no.~2, pp. 739--756, 2021.

\bibitem{D2}
K.~Baker, ``Learning warm-start points for {AC} optimal power flow,'' in
  \emph{IEEE 29th International Workshop on Machine Learning for Signal
  Processing (MLSP)}.\hskip 1em plus 0.5em minus 0.4em\relax IEEE, 2019, pp.
  1--6.

\bibitem{thomas}
T.~Falconer and L.~Mones, ``Leveraging power grid topology in machine learning
  assisted optimal power flow,'' \emph{IEEE Trans. Power Syst.}, vol.~38,
  no.~3, pp. 2234--2246, 2023.

\bibitem{KyriBaker}
K.~Baker, ``Solutions of {DC OPF} are never {AC} feasible,'' in
  \emph{Proceedings of the Twelfth ACM International Conference on Future
  Energy Systems}, 2021.

\bibitem{ejor}
F.~Zohrizadeh, C.~Josz, M.~Jin, R.~Madani, J.~Lavaei, and S.~Sojoudi, ``A
  survey on conic relaxations of optimal power flow problem,'' \emph{Eur. J.
  Oper. Res.}, vol. 287, no.~2, pp. 391--409, 2020.

\bibitem{Schapire}
D.~A. McAllester and R.~E. Schapire, ``On the convergence rate of good-turing
  estimators.'' in \emph{COLT}, 2000, pp. 1--6.

\bibitem{ML}
E.~Alpaydin, \emph{Introduction to machine learning}.\hskip 1em plus 0.5em
  minus 0.4em\relax MIT press, 2020.

\bibitem{kernel2}
T.~Hastie, R.~Tibshirani, J.~H. Friedman, and J.~H. Friedman, \emph{The
  elements of statistical learning: data mining, inference, and
  prediction}.\hskip 1em plus 0.5em minus 0.4em\relax Springer, 2009.

\bibitem{Soroush1}
S.~Shafieezadeh-Abadeh, D.~Kuhn, and P.~M. Esfahani, ``Regularization via mass
  transportation,'' \emph{Journal of Machine Learning Research}, vol.~20, no.
  103, pp. 1--68, 2019.

\bibitem{Soroush2}
D.~Kuhn, P.~M. Esfahani, V.~A. Nguyen, and S.~Shafieezadeh-Abadeh,
  ``Wasserstein distributionally robust optimization: Theory and applications
  in machine learning,'' in \emph{Operations Research \& Management Science in
  the Age of Analytics}, 2019, pp. 130--166.

\bibitem{patternreco}
C.~M. Bishop and N.~M. Nasrabadi, \emph{Pattern recognition and machine
  learning}.\hskip 1em plus 0.5em minus 0.4em\relax Springer, 2006, vol.~4,
  no.~4.

\bibitem{vapnik}
V.~N. Vapnik, \emph{The nature of statistical learning theory}.\hskip 1em plus
  0.5em minus 0.4em\relax Springer, 1999.

\bibitem{kernel1}
K.-R. Muller, S.~Mika, G.~Ratsch, K.~Tsuda, and B.~Scholkopf, ``An introduction
  to kernel-based learning algorithms,'' \emph{IEEE Trans. Neural Netw.},
  vol.~12, no.~2, pp. 181--201, 2001.

\bibitem{Keerthi}
S.~S. Keerthi and C.-J. Lin, ``Asymptotic behaviors of support vector machines
  with {G}aussian kernel,'' \emph{Neural Computation}, vol.~15, no.~7, pp.
  1667--1689, 2003.

\bibitem{GitHub}
\BIBentryALTinterwordspacing
 [Online]. Available: \url{https://github.com/farzanehpourahmadi/UCP.git}
\BIBentrySTDinterwordspacing

\bibitem{UC_AC2}
A.~Castillo, C.~Laird, C.~A. Silva-Monroy, J.-P. Watson, and R.~P. O’Neill,
  ``The unit commitment problem with ac optimal power flow constraints,''
  \emph{IEEE Trans. Power Syst.}, vol.~31, no.~6, pp. 4853--4866, 2016.

\bibitem{UC_AC3}
D.~Tuncer and B.~Kocuk, ``An {MISOCP}-based decomposition approach for the unit
  commitment problem with {AC} power flows,'' \emph{IEEE Trans. Power Syst.},
  vol.~38, no.~4, pp. 3388--3400, 2023.

\bibitem{UC_AC4}
V.~K. Tumuluru and D.~H.~K. Tsang, ``A two-stage approach for network
  constrained unit commitment problem with demand response,'' \emph{IEEE Trans.
  Smart Grid}, vol.~9, no.~2, pp. 1175--1183, 2018.

\bibitem{UC_AC5}
S.~Dehghan, P.~Aristidou, N.~Amjady, and A.~J. Conejo, ``A distributionally
  robust {AC} network-constrained unit commitment,'' \emph{IEEE Trans. Power
  Syst.}, vol.~36, no.~6, pp. 5258--5270, 2021.

\bibitem{jabr}
R.~A. Jabr, ``A conic quadratic format for the load flow equations of meshed
  networks,'' \emph{IEEE Trans. Power Syst.}, vol.~22, no.~4, pp. 2285--2286,
  2007.

\end{thebibliography}
